\documentclass{article}
\usepackage[utf8]{inputenc}

\usepackage{fullpage}
\usepackage{amsmath}
\usepackage{amssymb}
\usepackage{graphicx}
\usepackage{subcaption}
\usepackage{comment}
\usepackage{soul, xcolor}
\usepackage{authblk}
\usepackage{algorithm}
\usepackage{algorithmic}
\usepackage{diagbox}

\usepackage[colorlinks]{hyperref} 

\usepackage[
    backend=biber,
    style=numeric,
    sortcites=true,
    sorting=none,
    hyperref=true,
    natbib,
]{biblatex}

\usepackage{biblatex}

\hypersetup{
    citecolor = blue,
}

\graphicspath{{Figures/}}

\title{Scheduled Relaxation Jacobi schemes for non-elliptic partial differential equations}

\author[,a]{Mohammad Shafaet Islam\thanks{Corresponding author \vfill \ \ \textit{Email addresses}: \texttt{moislam@mit.edu} (Mohammad Shafaet Islam), \texttt{qiqi@mit.edu} (Qiqi Wang)}}
\author[a]{Qiqi Wang}

\affil[a]{\footnotesize{\textit{Department of Aeronautics and Astronautics, Massachusetts Institute of Technology, 77 Massachusetts Avenue, Cambridge, MA 02139, USA}}}

\date{}

\begin{document}

\maketitle

\begin{center}
\rule{\textwidth}{.5pt}  
\end{center}
\begin{abstract}

The Scheduled Relaxation Jacobi (SRJ) method is a linear solver algorithm which greatly improves the convergence of the Jacobi iteration through the use of judiciously chosen relaxation factors (an SRJ scheme) which attenuate the solution error. Until now, the method has primarily been used to accelerate the solution of elliptic PDEs (e.g. Laplace, Poisson's equation) as the currently available schemes are restricted to solving this class of problems. The goal of this paper is to present a methodology for constructing SRJ schemes which are suitable for solving non-elliptic PDEs (or equivalent, nonsymmetric linear systems arising from the discretization of these PDEs), thereby extending the applicability of this method to a broader class of problems. These schemes are obtained by numerically solving a constrained minimization problem which guarantees the solution error will not grow as long as the linear system has eigenvalues which lie in certain regions of the complex plane. We demonstrate that these schemes are able to accelerate the convergence of standard Jacobi iteration for the nonsymmetric linear systems arising from discretization of the 1D and 2D steady advection-diffusion equations.

\vspace{0.1in}

\end{abstract}
\begin{center}
\rule{\textwidth}{.5pt}  
\end{center}

\section{Introduction}
\label{sec:intro}

Scientists and engineers are often interested in the simulation of certain physical phenomena which are governed by partial differential equations (PDEs). These PDEs govern a wide variety of physical phenomena such as fluid flow \cite{cfd2030}, heat transfer \cite{jaluria2017computational} and electromagnetics \cite{electromagnetics-book}. Since these equations are generally not solvable analytically (except in the case of extremely idealized assumptions which may not be representative of real world settings), a numerical method is usually employed in order to solve these PDEs computationally. These numerical discretization methods lead to a large sparse system of equations that must be solved for the solution at specific points, cells or nodes in the domain \cite{Saad03}. However, the solution of these linear systems of equations can be a major bottleneck of this simulation procedure. The development of efficient linear solver methods can accelerate the study of physical phenomena of interest to scientists and engineers and improve our understanding of the natural world.

A plethora of linear solver algorithms exist for solving linear systems of equations \cite{forsythe}. The Jacobi iteration method \cite{jacobi} is perhaps the simplest linear solver in the class of iterative linear solvers. One key characteristic of the method is its inherently parallel nature. In the Jacobi update, each degree of freedom of the solution vector can be updated independently, making it an ideal candidate for implementation on high performance computing architectures (such as GPUs). For this reason, the method provides a good starting point towards the development of an efficient linear solver. One drawback of the Jacobi iteration is that it may be slow to converge, particularly for stiff problems (e.g. Poisson's equation on heavily refined grids). To ameliorate these convergence issues, Yang and Mittal developed the Scheduled Relaxation Jacobi (SRJ) method \cite{Yang14}. 

The Scheduled Relaxation Jacobi method aims to improve the convergence of the standard Jacobi iteration by applying a set of relaxed Jacobi updates with well chosen relaxation factors. The method is in the vein of polynomial acceleration methods to accelerate the convergence of stationary iteration methods \cite{varga1962iterative}. This includes the Chebyshev semi-iterative method \cite{golub1961chebyshev}, as well as work by Richardson \cite{richardson} and Young \cite{young1954iterative} to accelerate linear solver convergence behavior through the use of relaxation. At each cycle of the SRJ method, a fixed number of Jacobi iterations $M$ are performed with $P$ predetermined factors. These schemes aim to minimize the maximum achievable amplification for the 2D Laplace equation discretized on various grid sizes. In particular, Yang and Mittal derived schemes with $P = 2,3,4,5$ relaxation factors for the 2D Laplace equation discretized with $N = 16, 32, 64, 128, 256, 512$ degrees of freedom in each dimension. They observe speedup factors upwards of 100 times from the use of certain schemes relative to Jacobi iteration. Further work on the SRJ method was conducted by Adsuara et al. They were able to analytically simplify the minimization problem considered by Yang so that schemes associated with larger $P$ (up to 15) for much larger grid resolutions ($N$  up to $2^{15}$) can be derived \cite{Adsuara1}. Later, Adsuara et al. proposed further improvements to the method, such as specifying that each relaxation factor is used only once in a given cycle so that $P = M$. Under this new constraint, they are able to develop new schemes which demonstrate improved convergence when solving the 2D Laplace equation, compared to the original schemes derived earlier. These new schemes are also easier to obtain as the relaxation factors can be found given the roots of the Chebyshev polynomials, which are known analytically \cite{Adsuara2}. Islam and Wang extend this idea further, deriving SRJ schemes based on achieving solution error amplification related to the Chebyshev polynomials. They derive schemes which are generally suited to solving Poisson problems regardless of the specific discretization used, and employ a data driven approach to select which scheme to use at each cycle of the SRJ method during the linear solve procedure \cite{islam2020}. Many avenues for further augmenting the practical performance of these SRJ schemes are also being explored. These include integration of these schemes within multigrid solvers \cite{Yang17}, GPU implementation \cite{Adsuara3}, as well as other theoretical \cite{aasrj2019, pratapa2016anderson} and data-driven approaches \cite{islam2020, xiang2021neuroevolution} to improve their performance.

The SRJ method is a promising candidate as a high performance linear solver method for simulation. These schemes have demonstrated faster convergence relative to the Jacobi iteration, and with further augmentations from theory, high performance computing architectures, and data driven approaches can become an extremely powerful linear solver method. However, a current limitation is that the schemes that have currently been developed until now are restricted to solving PDEs of the elliptic type (e.g. Laplace, Poisson's equation). Although non-elliptic PDEs are also prevalent in nature (e.g. advection-dominated flows which are represented by the parabolic advection-diffusion equation), the current schemes are likely ill-suited for such problems. The goal of this paper is to develop SRJ schemes which would be appropriate for solving non-elliptic PDEs, or equivalently, the nonsymmetric linear systems arising from discretization of these PDEs. This broadens the applicability of the method to a wider class of problems which are of interest to scientists and engineers.


This paper is structured as follows. Section 2 presents a review of the derivation of SRJ schemes for solving symmetric linear systems arising from the discretization of elliptic PDEs. The analysis follows the schemes developed in \cite{islam2020}. Section 3 develops the methodology for deriving schemes which are appropriate for solving non-elliptic PDEs. The schemes are derived by minimizing the solution error amplification in elliptical regions of the complex plane, and is akin to previous effort to extend polynomial acceleration methods to nonsymmetric matrices \cite{manteuffel1977tchebychev}. Section 4 presents results on the convergence behavior of these new schemes for solving the one-dimensional and two-dimensional steady advection-diffusion equation which fall into the class of parabolic PDEs. These schemes demonstrate accelerated convergence compared to the standard Jacobi iteration, and continue to provide acceleration in cases where the original schemes developed solely for elliptic PDEs fail to converge. Section 5 provides concluding remarks for this work. 

\section{Derivation of SRJ schemes for symmetric linear systems}
\label{sec:background-on-srj-schemes}

This section provides background on the Scheduled Relaxation Jacobi schemes (sets of relaxation factors) which have been developed for solving elliptic PDEs (or equivalently, the symmetric linear systems arising from their discretization). The analysis here is based on the work in \cite{islam2020}. Further information regarding these schemes can be found in this reference.

Our goal is to solve a linear system of equations $Ax = b$, $A \in \mathbb{R}^{n \times n}$, $x \in \mathbb{R}^{n}$, $b \in \mathbb{R}^{n}$. Our starting point is the Jacobi iterative method, which is given by the update Equation \eqref{eqn:jacobi}. In this equation, $x^{(n+1)}$ is the solution at the next iteration, $x^{(n)}$ is the current solution to the linear system, and $L$, $U$, and $D$ are the lower triangular, upper triangular, and diagonal matrices, respectively, containing the corresponding elements of $A$. One can write down an iteration matrix associated with the iterative update denoted by $B_{\text{J}} = -D^{-1} (L+U)$.
\begin{equation}
    x^{(n+1)} =  \underbrace{-D^{-1} (L+U)}_{B_{\text{J}}}  x^{(n)} + D^{-1} b
    \label{eqn:jacobi}
\end{equation}
The Jacobi iterative method is known to converge for a linear system if the spectral radius of the iteration matrix is less than 1. This is expressed in Equation \eqref{eqn:spectral-radius}. The spectral radius also gives a good indication of the asymptotic convergence rate, with a lower spectral radius corresponding to faster convergence.
\begin{equation}
    \rho (B_{\text{J}}) < 1
    \label{eqn:spectral-radius}
\end{equation}
We can also define a relaxed Jacobi update, where the solution at the next step is equal to a weighted average of the solution value from the Jacobi update weighted by some factor $\omega$, and the current solution weighted by $(1-\omega)$. This relaxed Jacobi update is given by Equation \eqref{eqn:relaxed-jacobi}.
\begin{equation}
    x^{(n+1)} = \left[ (1-\omega) I + \omega  B_{\text{J}} \right] x^{(n)} + \omega D^{-1} b 
    \label{eqn:relaxed-jacobi}
\end{equation}
The relaxation factor employed affects the convergence of the Jacobi iteration. To see this, we first define the solution error vector at step $n$ as $e^{(n)} \equiv x^{(n)} - x$ where $x$ is the exact solution to the linear system. The exact solution satisfies the update equation \eqref{eqn:relaxed-jacobi} exactly. Given this, we may derive an update equation for the evolution of the error from one step to the next. This results in the error evolution equation at each iteration given by Equation \eqref{eqn:error-update}. The accumulation of the solution error is based on the matrix $B_{\omega} = (1-\omega) I + \omega  B_{\text{J}}$. Convergence of weighted Jacobi is guaranteed if the matrix $B_{\omega}$ has a spectral radius less than 1 (i.e. $\rho(B_{\omega}) < 1$). Furthermore, the asymptotic convergence of the relaxed Jacobi update is related to this spectral radius. The relaxation factor impacts the spectral radius of the matrix $B_{\omega}$, thereby affecting the convergence rate of the relaxed Jacobi iterations. In particular, it is beneficial to apply a relaxation factor which reduces the spectral radius of $B_{\omega}$, with the caveat that applying overrelaxation alone ($\omega > 1$) is known to cause Jacobi iteration to diverge.
\begin{equation}
    e^{(n+1)} = \underbrace{\left[ (1-\omega) I + \omega  B_{\text{J}} \right]}_{B_{\omega}} e^{(n)} 
    \label{eqn:error-update}
\end{equation}
We now consider an iteration scheme where $M$ iterations of the relaxed Jacobi method are performed, each with a distinct relaxation factor denoted by $\omega_i$. Let these iterations comprise one cycle of the SRJ algorithm, and denote the overall iteration matrix associated with these $M$ iterations by  $B_{\text{SRJ}}$. We also denote the iteration matrix associated with each individual iteration by $B_{\omega_{i}}$. If $e^{(n)}$ and $e^{(n+1)}$ now represent the error prior to and after a cycle of $M$ iterations, then the solution error evolves as follows
\begin{equation}
    e^{(n+1)} = B_{\text{SRJ}} e^{(n)}  = \prod_{i=1}^{M} B_{\omega_i} e^{(n)}  = \prod_{i=1}^{M} \left[ (1-\omega_i) I + \omega_i  B_{\text{J}} \right] e^{(n)}
\end{equation}
In order for the error to decay from one SRJ cycle to the next, the spectral radius of the SRJ iteration matrix must be less than 1 (i.e. $\rho(B_{\text{SRJ}}) < 1$). The error analysis in Equations \eqref{eqn:error-analysis-1} - \eqref{eqn:error-analysis-end} illustrates this. Here, we express the initial error vector $e^{(0)}$ in a basis comprised of the eigenvectors of the matrix $B_{\text{SRJ}}$, which are denoted by $v_{j}$. We also use the fact that if $\lambda_{\text{SRJ}}$ are the corresponding eigenvalues of $B_{\text{SRJ}}$, then $B_{\text{SRJ}} v_{j}= \lambda_{\text{SRJ}} v_{j}$.
\begin{align}
    \label{eqn:error-analysis-1}
    e^{(n+1)} & = B_{\text{SRJ}} e^{(n)}  \\ 
                      &= \left( B_{\text{SRJ}} \right)^{n+1} e^{(0)}  \\ 
                      &= \left( B_{\text{SRJ}} \right)^{n+1} \sum_{j} a_{j} v_{j}  \\ 
                      &= \sum_{j} a_{j} \left( B_{\text{SRJ}} \right)^{n+1} v_{j}  \\
                      &= \sum_{j} a_{j} (\lambda_{\text{SRJ}})^{n+1} v_{j}  
    \label{eqn:error-analysis-end}
\end{align}
The error analysis above implies that the error will grow unbounded if the magnitude of any of the eigenvalues $\lambda_{\text{SRJ}}$ are greater than 1. The relaxation factors for the Scheduled Relaxation procedure must therefore be chosen to ensure that the resulting eigenvalues have magnitude less than 1 (i.e. the spectral radius of the SRJ iteration matrix is less than 1). We can show that the eigenvalues of the SRJ iteration matrix (denoted by $\lambda_{\text{SRJ}}$) are related to the eigenvalues of the Jacobi iteration matrix (which we will denote by $\lambda_{\text{J}}$) as follows. 
If $v_{j}$ is an eigenvector of both $B_{\text{SRJ}}$ and $B_{\text{J}}$, then
\begin{align}
    B_{\text{SRJ}} v_{j} &= \prod_{i = 1}^{M} \left[ (1-\omega_{i}) I + \omega_{i} B_{\text{J}} \right] v_{j} \\ 
    &= \prod_{i = 1}^{M} \left[ (1-\omega_{i}) v_{j} + \omega_{i} \lambda_{\text{J}} v_{j} \right]  \\ 
    &= \prod_{i = 1}^{M} \left[ (1-\omega_{i}) + \omega_{i} \lambda_{\text{J}} \right] v_{j} 
    \label{eqn:eigen-srj-1}
\end{align}
Additionally, by definition
\begin{equation}
     B_{\text{SRJ}} v_{j} = \lambda_{\text{SRJ}} v_{j}
     \label{eqn:eigen-srj-2}
\end{equation}
Given Equations \eqref{eqn:eigen-srj-1} and \eqref{eqn:eigen-srj-2}, we obtain the following relationship between the eigenvalues of the SRJ iteration matrix and those of the Jacobi iteration matrix
\begin{equation}
    \lambda_{\text{SRJ}} = G_M(\lambda_{\text{J}})\;,\quad\mbox{where}\quad G_M(\lambda):= \prod_{i=1}^{M} \left[(1-\omega_i) + \omega_i \lambda\right]
    \label{eqn:amplification-M-scheme}
\end{equation}
According to Equation \eqref{eqn:amplification-M-scheme}, the eigenvalues $\lambda_{\text{SRJ}}$  are an $M$-degree polynomial of the eigenvalues $\lambda_{\text{J}}$, which we will call the amplification polynomial and denote by $G_{M}$. The Jacobi iteration will converge as long as all eigenvalues $\lambda_{\text{J}}$ lie in $(-1,1)$. For the SRJ iterations to converge under the same conditions, the magnitude of the amplification polynomial $G_{M}(\lambda)$ must be bounded by 1 for all inputs $\lambda \in (-1,1)$. 

We can also determine the asymptotic convergence rate of the SRJ method. The convergence rate of the Jacobi iteration is determined by the spectral radius of $B_{\text{J}}$, while the convergence rate of the SRJ method is determined by the spectral radius of $B_{\text{SRJ}}$. The spectral radius of the SRJ iteration matrix is found by evaluating the amplification polynomial $G_{M}$ at all of the Jacobi iteration matrix eigenvalues $\lambda_{\text{J}}$ (resulting in the SRJ iteration matrix eigenvalues $\lambda_{\text{SRJ}}$), and taking the largest absolute value of these eigenvalues. This is expressed in Equation \eqref{eqn:srj-convergence}. 
\begin{equation}
    \rho(B_{\text{SRJ}}) = \max \left|\lambda_{\text{SRJ}} \right| =  \max \left| G_{M}(\lambda_{\text{J}}) \right| 
    \label{eqn:srj-convergence}
\end{equation}
The analysis above suggests that the amplification polynomial should satisfy a few requirements for the associated SRJ scheme to converge. First of all, the amplification polynomial must be bounded between -1 and 1 for any input $\lambda \in (-1,1)$. Second, to ensure that the SRJ scheme converges efficiently, the spectral radius $\rho(B_{\text{SRJ}})$ should be as small as possible. This motivates the construction of amplification polynomials which are bounded for as large a range of $\lambda \in (-1,1)$ as possible. Given these requirements, we can construct a set of amplification polynomials for different orders $M$. Given an amplification polynomial, we can then derive the corresponding relaxation factors $\omega_{i}$ which yield this amplification behavior.

\begin{figure}[htbp!]
    \centering
    \includegraphics[width=0.6\textwidth]{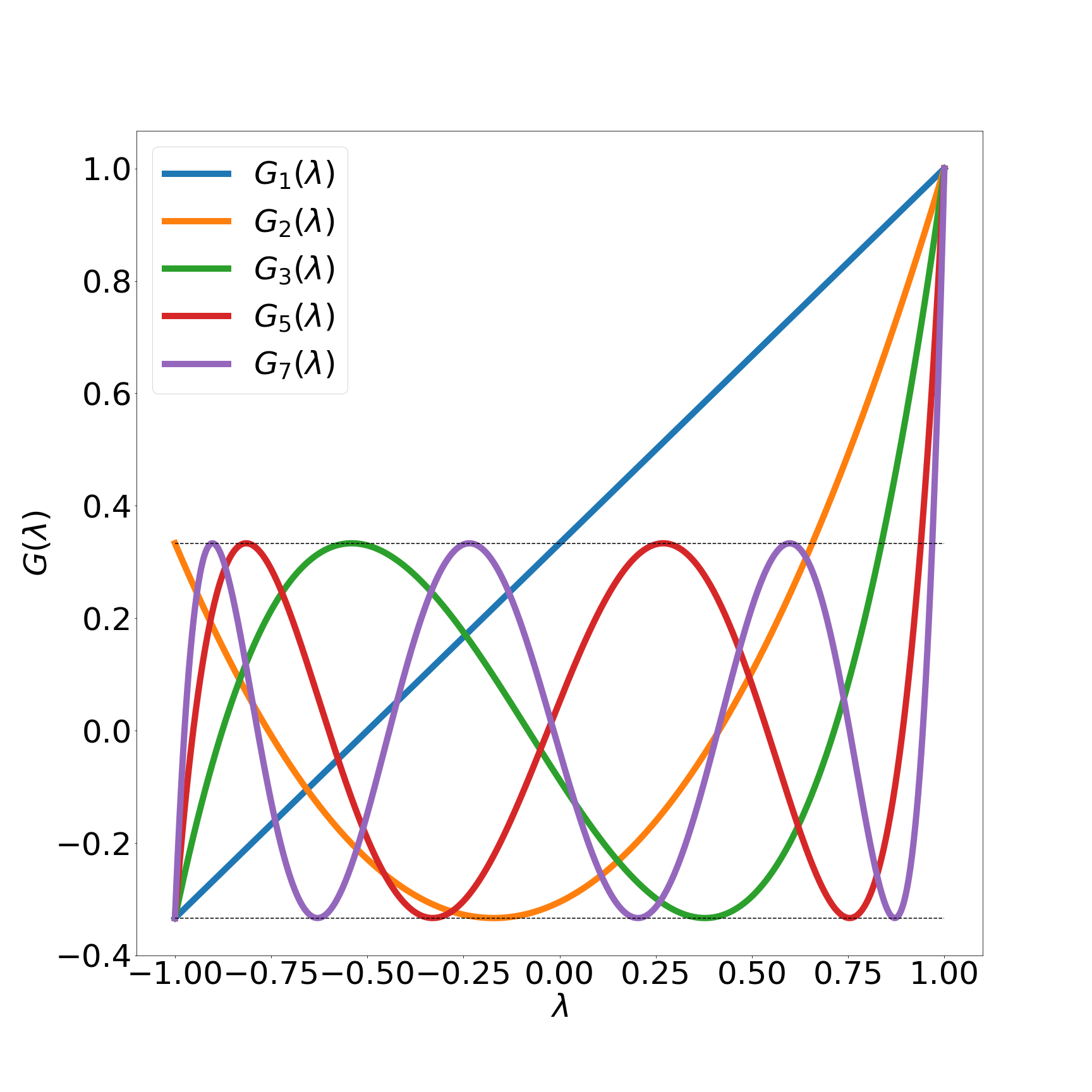}
    \caption{Amplification polynomials $G_{M}(\lambda)$ for $M = 1,2,3,5,7$. The polynomials are bounded by $\pm \frac{1}{3}$ for some region within $\lambda \in (-1,1)$, which grows as $M$ increases.}
    \label{fig:amplification_polynomials}
\end{figure}

Figure \ref{fig:amplification_polynomials} illustrates a few of the amplification polynomials which were designed for increasing order $M$, appropriate for problems of different stiffness. These amplification polynomials are the Chebyshev polynomials under a certain affine transformation which ensures that the amplification is bounded by a value of $\pm \frac{1}{3}$ for the widest range possible in $\lambda \in (-1,1)$. Specifically, given the $M$ order Chebyshev polynomial (denoted by $T_{M}(\lambda)$) the $M$ order amplification polynomial $G_{M}$ is found by 
\begin{equation}
    G_{M}(\lambda) = \frac{T_{M} (f(\lambda)) }{3} \;,\quad\mbox{where}\quad f(\lambda):= \frac{(\lambda^* + 1)\lambda + (\lambda^*-1)}{2}
    \label{eqn:transformation}
\end{equation}
where $\lambda^*$ satisfies $T_M (\lambda^*) = 3$. Table \ref{tab:amplification-polynomials-table} shows the functional representations of the amplification polynomials for $M = 1,2,3,5$ plotted in Figure \ref{fig:amplification_polynomials} as well as the corresponding Chebyshev polynomials. The maximum input $\lambda_{\text{max}}$ for which the amplification polynomial is bounded by $\frac{1}{3}$ is also shown (and grows closer to 1 with the polynomial order $M$). This $\lambda_{\text{max}}$ can be thought of as the input 1 in the Chebyshev polynomial (the upper bound of the inputs where the Chebyshev polynomial is bounded) mapped to a corresponding input for the amplification polynomial under the linear transformation. This can be written as
\begin{equation}
    \lambda_{\text{max}} = g(1) \;,\quad\mbox{where}\quad g(\lambda) := f^{-1}(\lambda) = \frac{2 \lambda}{1+\lambda^{*}} + \frac{1 - \lambda^{*}}{1+\lambda^{*}}
    \label{eqn:transformation}
\end{equation}
Based on this, $\lambda_{\text{max}}$ is given by
\begin{equation}
    \lambda_{\text{max}} = \frac{3-\lambda^{*}}{1+\lambda^{*}}
    \label{eqn:transformation}
\end{equation}
As $M$ increases, the $\lambda^{*}$ such that $T_{M}(\lambda^*) = 3$ gets closer to 1, so $\lambda_{\text{max}}$ also grows closer to 1 (as implied by Table \ref{tab:amplification-polynomials-table}). This shows that increasing $M$ increases the overall range of $\lambda \in (-1,1)$ for which the amplification is bounded by $\frac{1}{3}$. 

\begin{table}[htbp!]
    \caption{The amplification polynomials $G_{M}$ and the maximum $\lambda$ for which the polynomials are bounded by $\frac{1}{3}$, for $M = 1,2,3,5$. The corresponding Chebyshev polynomials $T_{M}(\lambda)$ are also shown.}
    \centering
    \begin{tabular}{|c|l|l|c|}
        \hline
        $M$ & $T_{M}(\lambda)$ & $G_{M}(\lambda)$ & $\lambda_{\text{max}}$  \\
        \hline
        1 & $\lambda$ & $\frac{2}{3}\lambda + \frac{1}{3}$ & 0.0 \\
        \hline
        2 & $2\lambda^2-1$ & $\frac{2}{3} (1.2071 \lambda + 0.2071)^2 - \frac{1}{3}$ & 0.6569 \\
        \hline
        3 & $4\lambda^3-3\lambda$ & $\frac{4}{3} (1.0888 \lambda + 0.0888)^3 - (1.0888 \lambda + 0.0888)$ & 0.8368 \\
        \hline
        5 & $16\lambda^5 - 20\lambda^3 + 5\lambda$ & $\frac{16}{3}(1.0314 \lambda + 0.0314)^5 - \frac{20}{3}(1.0314 \lambda + 0.0314)^3 + \frac{5}{3}(1.0314 \lambda + 0.0314)$ & 0.9391 \\
        \hline
    \end{tabular}
    \label{tab:amplification-polynomials-table}
\end{table}

The SRJ scheme parameters corresponding to the amplification polynomials shown in Figure \ref{fig:amplification_polynomials} are listed in Table \ref{tab:srj-schemes}. Reference \cite{islam2020} provides more information on obtaining the affine transformation from Chebyshev to amplification polynomials, as well as obtaining the SRJ scheme parameters given the amplification polynomials. 

\begin{table}[htbp!]
    \caption{SRJ Relaxation factors associated with amplification for $M = 1,2,3,5,7$, shown in Figure \ref{fig:amplification_polynomials}}
    \centering
    \begin{tabular}{|c|l|}
        \hline
        $M$ & SRJ scheme parameters \\
        \hline
        1 & 0.66666667 \\
        \hline
        2 & 1.70710678, 0.56903559 \\
        \hline
        3 & 3.49402108, 0.53277784, 0.92457411 \\
        \hline
        5 & 9.23070105, 0.51215173, 0.97045899, 0.62486988, 2.1713295 \\
        \hline
        7 & 17.84007924, 0.50624677, 0.9845549, 1.69891732, 0.56014439, 4.06304526, 0.69311375\\
        \hline
    \end{tabular}
    \label{tab:srj-schemes}
\end{table}

\begin{figure}[htbp!]
    \begin{subfigure}{.5\textwidth}
        \centering
        \includegraphics[width=0.9\textwidth]{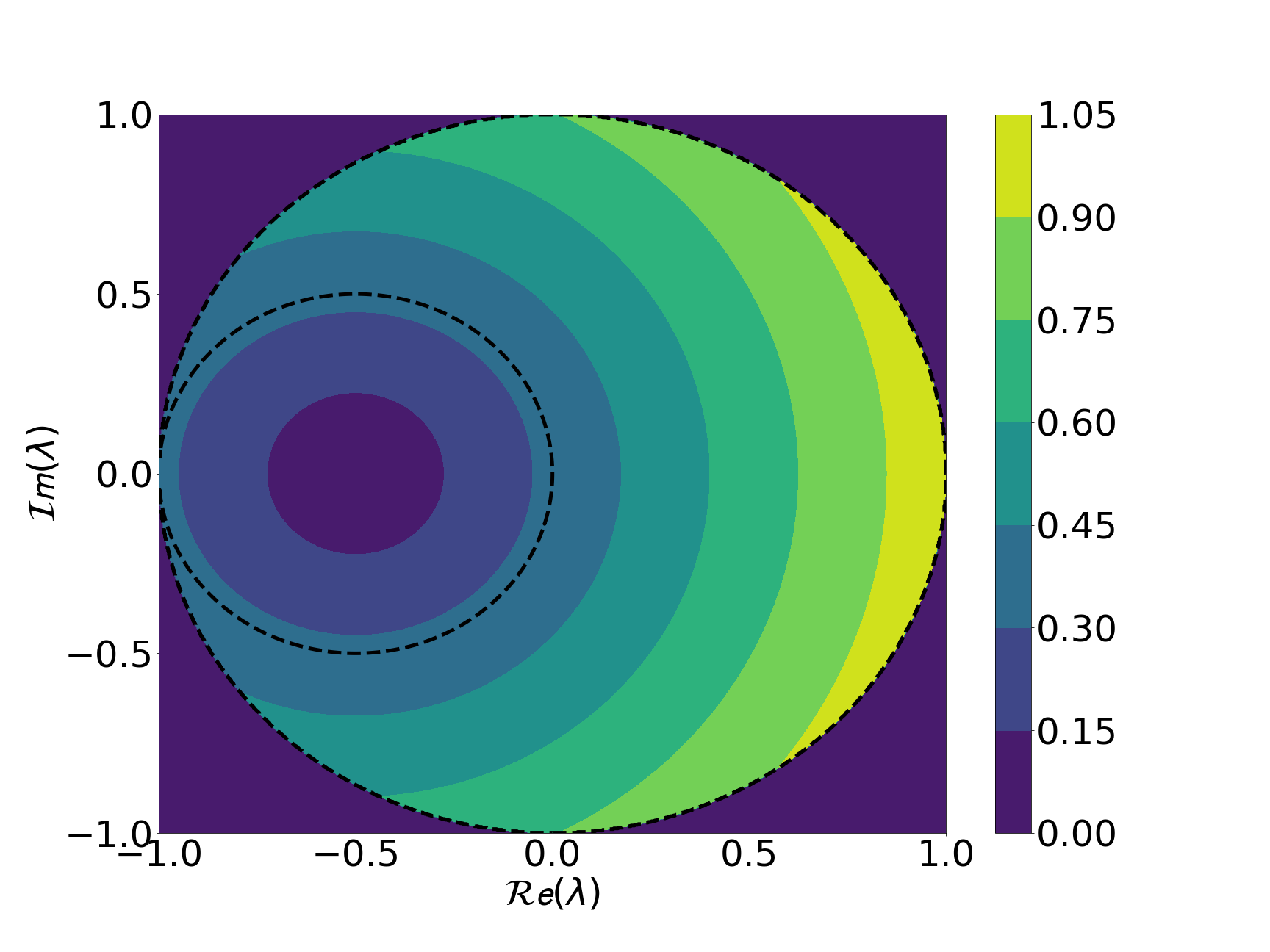}
        \caption{$M = 1$}
        \label{fig:M1_complex}
    \end{subfigure}
    \begin{subfigure}{.5\textwidth}
        \centering
        \includegraphics[width=0.9\textwidth]{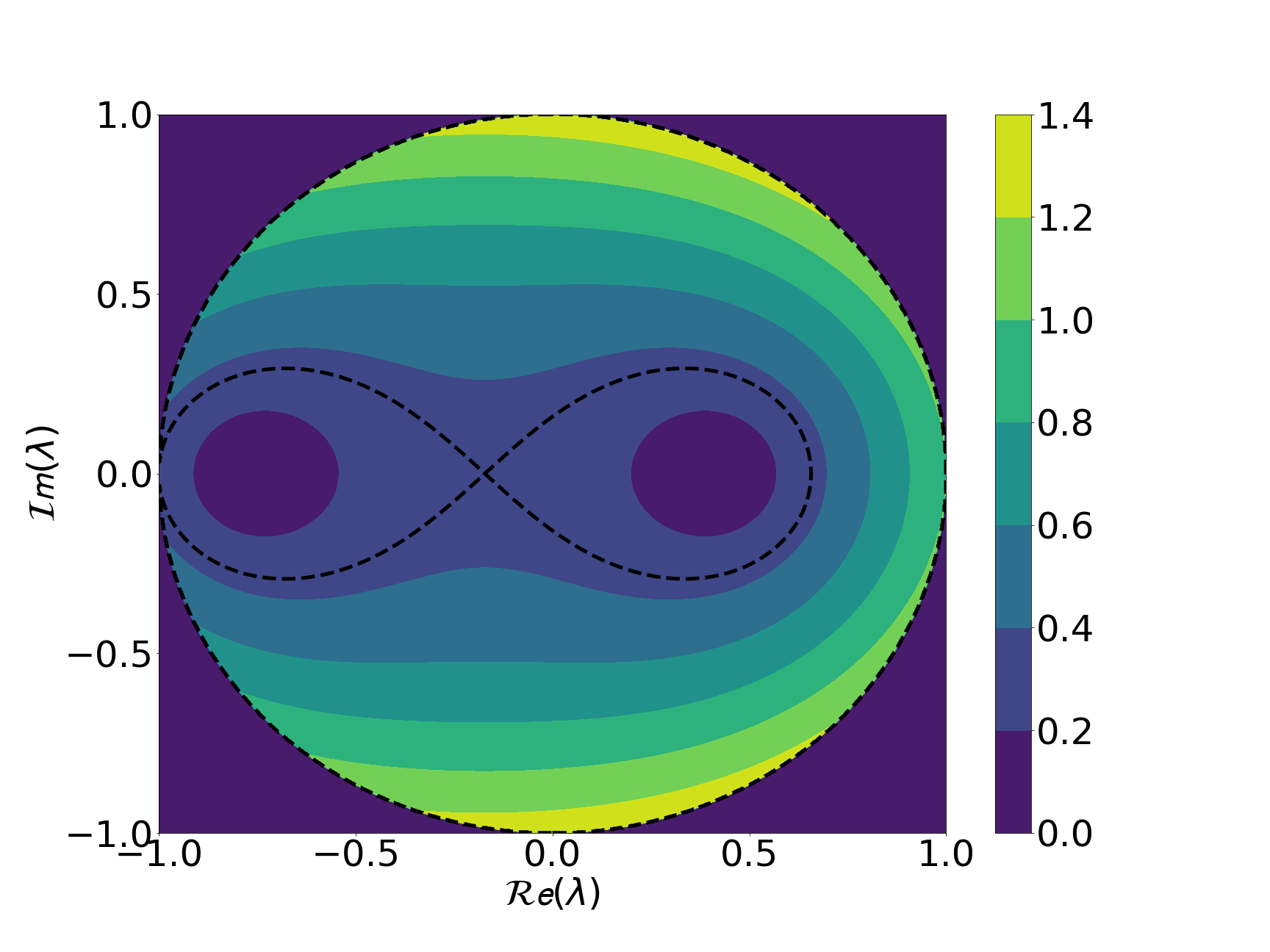}
        \caption{$M = 2$}
        \label{fig:M2_complex}
    \end{subfigure}
    \\
    \begin{subfigure}{.5\textwidth}
        \centering
        \includegraphics[width=0.9\textwidth]{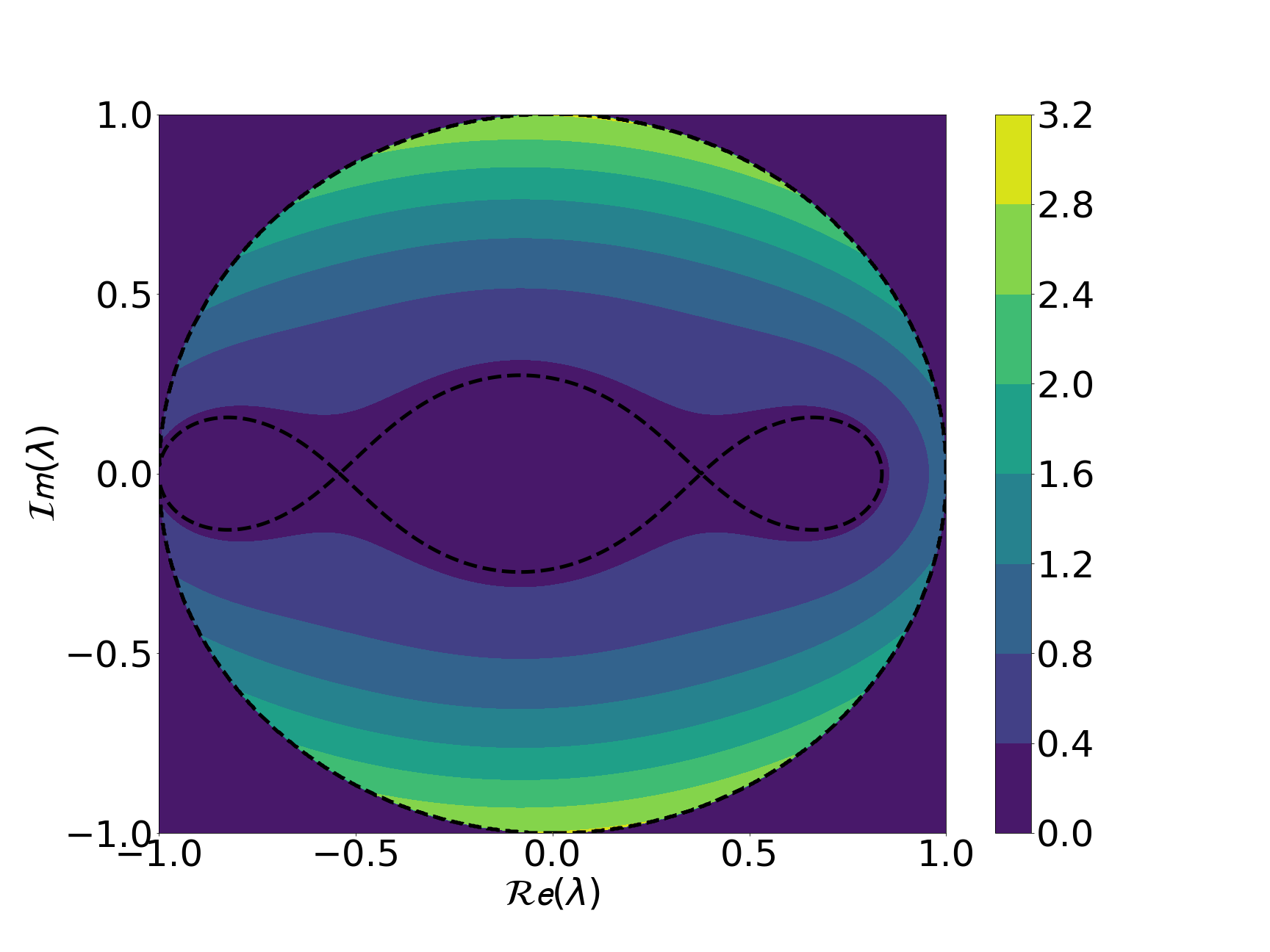}
        \caption{$M = 3$}
        \label{fig:M3_complex}
    \end{subfigure}
    \begin{subfigure}{.5\textwidth}
        \centering
        \includegraphics[width=0.9\textwidth]{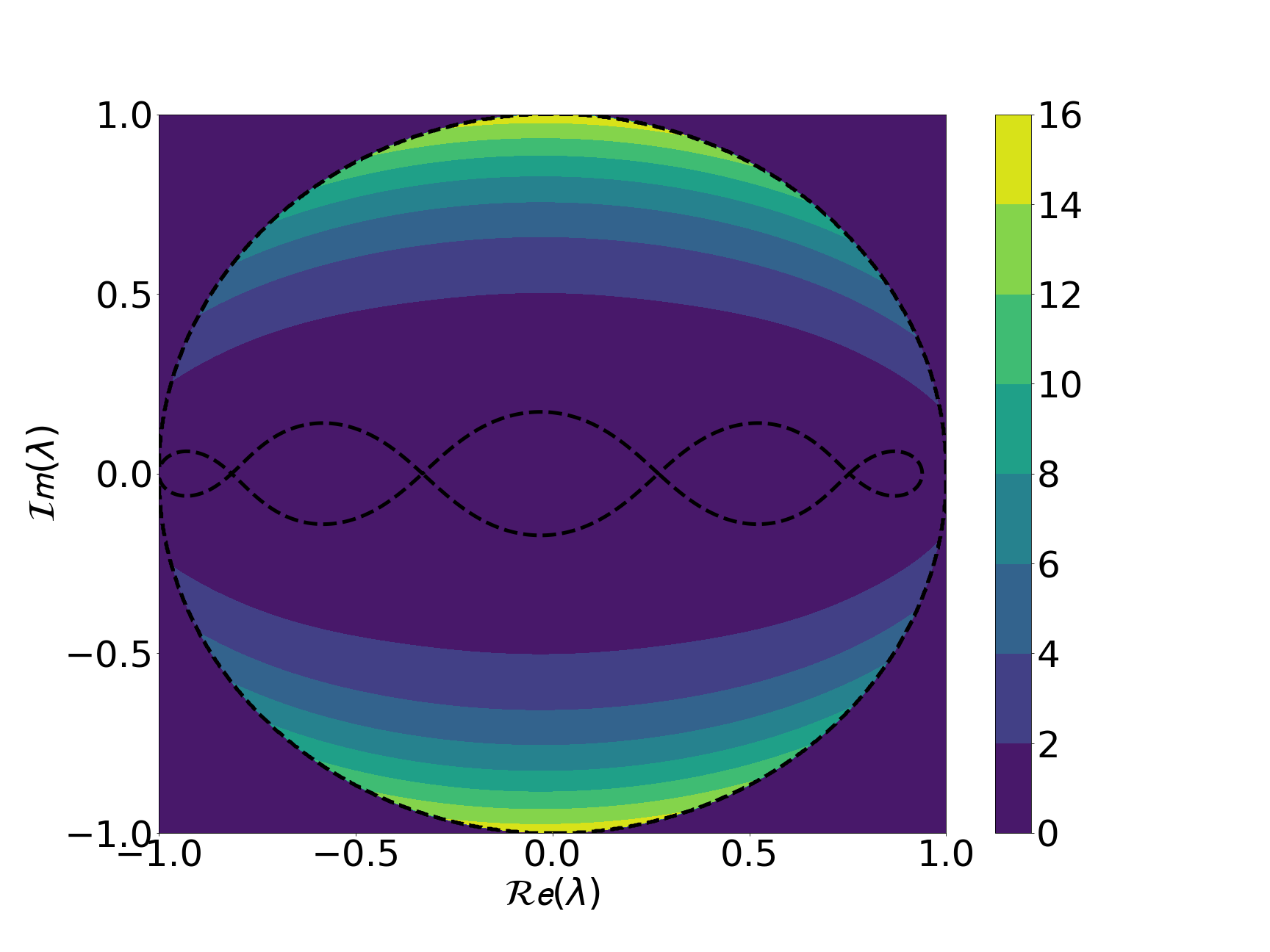}
        \caption{$M = 5$}
        \label{fig:M5_complex}
    \end{subfigure}
    \\
    \begin{center}
        \begin{subfigure}{.5\textwidth}
            \centering
            \includegraphics[width=0.9\textwidth]{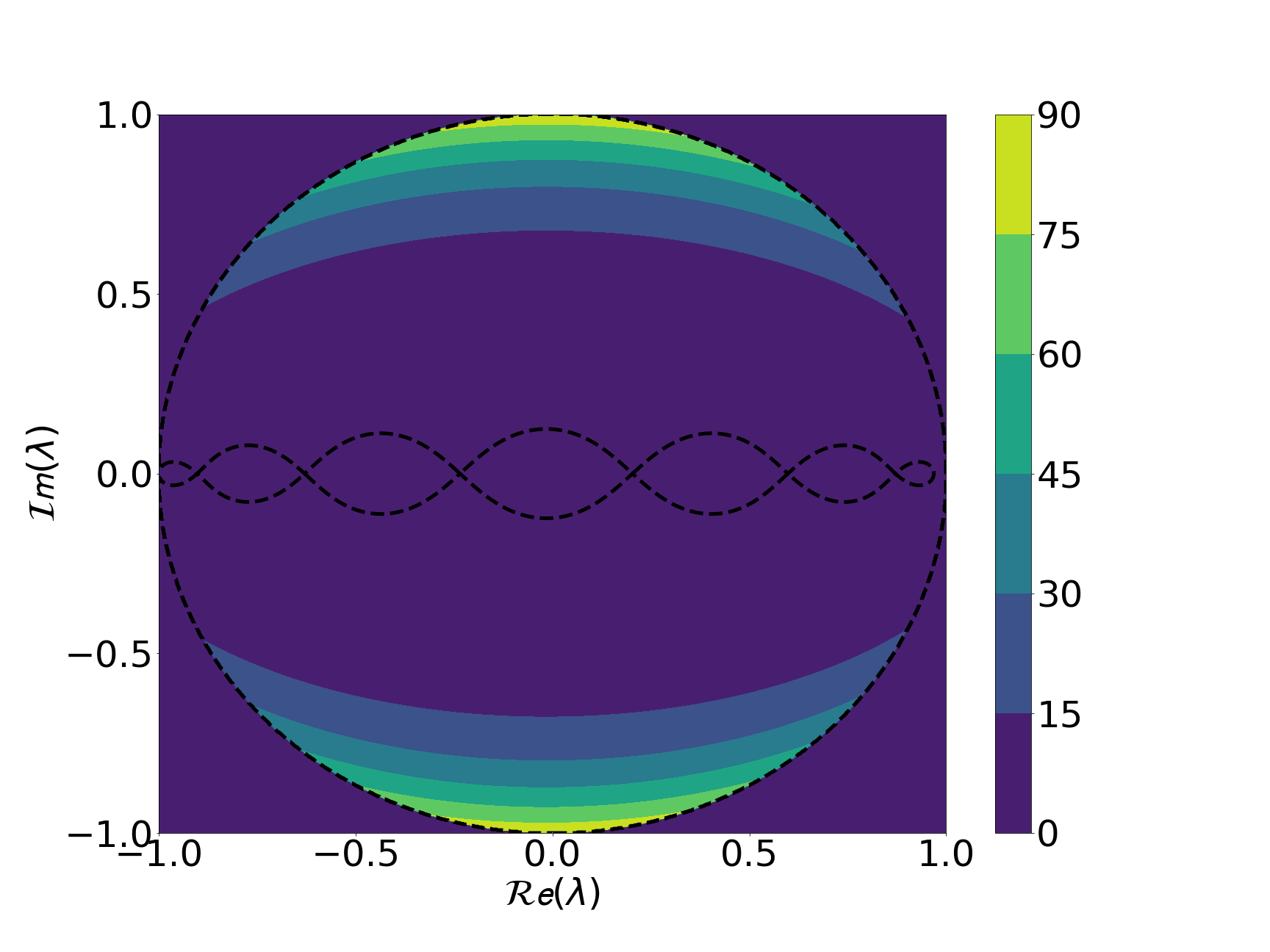}
            \caption{$M = 7$}
            \label{fig:M7_complex}
        \end{subfigure}
    \end{center}
    \caption{Contour plots of the amplification polynomials shown in Figure \ref{fig:amplification_polynomials} in the unit circle in the complex plane. The black dashes lines indicate contours where the polynomial is bounded by $\frac{1}{3}$. The amplification polynomials are not bounded by 1 within the unit circle. Therefore, given a linear system whose Jacobi iteration matrix eigenvalues lie anywhere in the unit circle, it is possible for the SRJ schemes to diverge even when the Jacobi iteration converges. As $M$ increases, the amplification magnitude exceeds 1 for a larger portion of the unit circle so the schemes are more likely to diverge. This motivates the need for new schemes for nonsymmetric matrices.}
    \label{fig:amplification_comparison}
\end{figure}

The schemes shown in Table \ref{tab:srj-schemes} are expected to attenuate the solution error as long as the Jacobi iteration matrix eigenvalues lie in $(-1,1)$ on the real axis. However, the eigenvalues of the iteration matrix lie on the real axis only when the matrix $A$ is symmetric. If $A$ is nonsymmetric, the iteration matrix may have complex eigenvalues. In this case, the Jacobi iteration will still converge if these iteration matrix eigenvalues lie within the unit circle in the complex plane. However, it is not guaranteed that our SRJ schemes can converge in such cases. To guarantee that our SRJ schemes will converge whenever Jacobi iteration converges, the associated amplification polynomials must not only be bounded on the real axis from $(-1,1)$, but more generally be bounded within the unit circle of the complex plane. Figures \ref{fig:M1_complex}-\ref{fig:M7_complex} illustrate contour plots of the magnitude of the amplification polynomials in Figure \ref{fig:amplification_polynomials} in the unit circle in the complex plane. Dashed lines in the figures represent an additional contour where the amplification polynomial is bounded by our bounding value of $\frac{1}{3}$.

Figures \ref{fig:M1_complex}-\ref{fig:M7_complex} illustrate that the magnitude of the amplification polynomials are not bounded by 1 in the unit circle. As the polynomial order increases, the maximum amplification within the domain is much higher, and the region in which the amplification is bounded by 1 is smaller. Additionally, the region where the polynomial is bounded by $\frac{1}{3}$ gets progressively closer to the real line with increasing $M$. This indicates that schemes associated with larger $M$ are even worse at handling complex eigenvalues compared to schemes with smaller $M$. The practical effect of the amplification polynomial exceeding 1 in the unit circle is the following. Suppose we wish to solve a linear system using a particular SRJ scheme. If any of the Jacobi iteration matrix eigenvalues lie in regions of the complex unit circle where the scheme's amplification magnitude exceeds 1, the SRJ scheme would most likely diverge while the Jacobi iteration would converge. Therefore, the approach to deriving SRJ schemes discussed above may be inadequate for solving nonsymmetric matrices whose eigenvalues are complex. This motivates the development of schemes with amplification behavior which is bounded in the complex plane, and therefore suitable for solving nonsymmetric linear systems of equations.

\section{Extension of SRJ schemes to nonsymmetric linear systems}
\label{sec:srj-for-nonsymmetric-schemes}

The goal of this section is to develop a methodology for deriving Scheduled Relaxation Jacobi schemes which are appropriate for solving nonsymmetric linear systems. The original schemes presented in Section \ref{sec:background-on-srj-schemes} are not guaranteed to converge when the eigenvalues of the Jacobi iteration matrix are complex (true for nonsymmetric linear systems). The approach presented here allows us to derive a family of schemes which are appropriate for solving linear systems whose Jacobi iteration matrix eigenvalues have different distributions in the complex plane.

\subsection{SRJ schemes as solutions to an optimization problem}

As a step towards obtaining Scheduled Relaxation Jacobi schemes for general nonsymmetric matrices, we first review the approach taken to derive the SRJ schemes for the symmetric matrices and discuss the key features of these schemes. The Scheduled Relaxation Jacobi schemes derived for the symmetric case were obtained by constructing a class of amplification polynomials which are able to attenuate the solution error. The main characteristic of the constructed amplification polynomials is that they are bounded by some bounding value (in our case $\pm \frac{1}{3}$) for the maximum possible range of eigenvalues within $\lambda \in (-1,1)$. For some general amplification, this range is $\lambda \in (-1,\lambda_{\text{max}})$ where $\lambda_{\text{max}}$ gets progressively close to 1 as the amplification polynomial order $M$ increases (see Table \ref{tab:amplification-polynomials-table}). Given this polynomial, the $M$ relaxation factors corresponding to this amplification can be obtained. 

Another viewpoint of the relaxation factors is that they are the scheme parameters which minimize the maximum value of the amplification polynomial within some bounds $\lambda \in (-1,\lambda_{\text{max}})$ where $\lambda_{\text{max}}$ is known for a given $M$ (according to Equation \eqref{eqn:transformation}). As a result, the scheme parameters can be posed as the solution to an optimization problem which seeks to minimize the maximum amplification within some region of the real axis, as written in equation \eqref{eqn:minimization-srj-real}. In equation \eqref{eqn:minimization-srj-real} $\mathcal{R}$ corresponds to the portion of the real axis $(-1, \lambda_{\text{max}})$ on which we bound the amplification polynomial $G_{M}$. For a given $M$, the solution to this optimization problem is the $M$ parameter SRJ scheme derived for symmetric matrices.
\begin{equation}
    \min_{\omega_i} \max_{\lambda \in \mathcal{R}} G_{M} (\lambda; \omega_i)  \;,\quad\mbox{where}\quad G_M(\lambda;\omega_i):= \prod_{i=1}^{M} \left[(1-\omega_i) + \omega_i \lambda\right]
    \label{eqn:minimization-srj-real}
\end{equation}

\begin{figure}[htbp!]
    \begin{subfigure}{.5\textwidth}
        \centering
        \includegraphics[width=\linewidth]{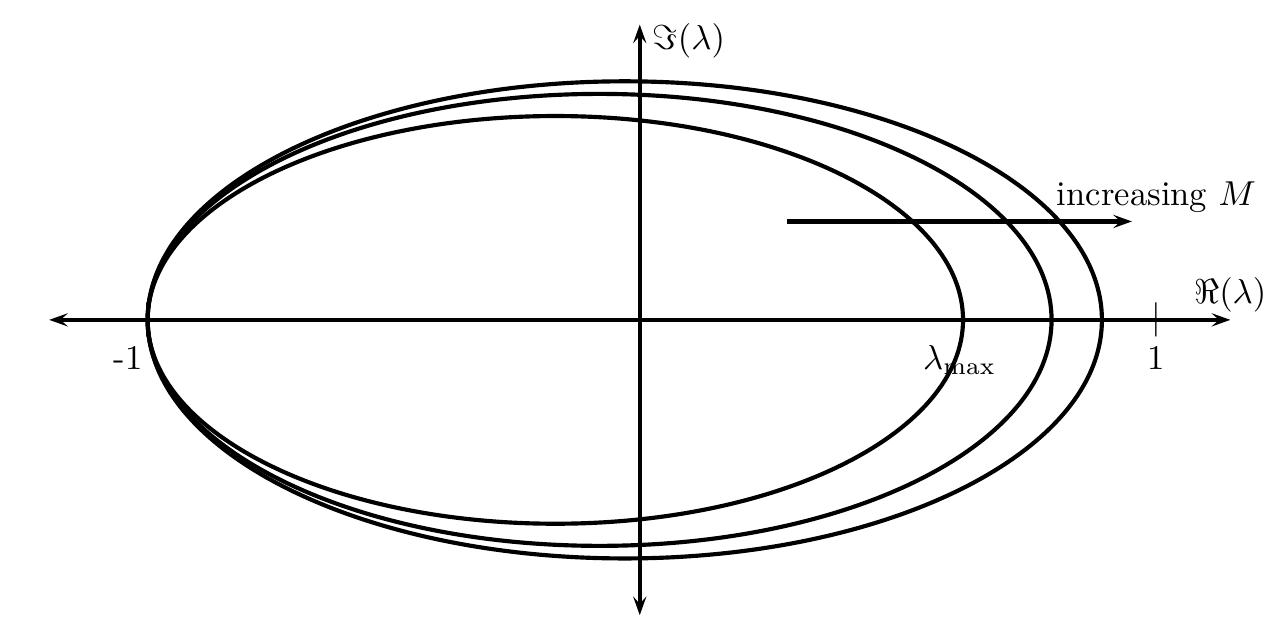}  
        \caption{Varying $M$}
        \label{fig:vary M}
    \end{subfigure}
    \begin{subfigure}{.5\textwidth}
        \centering
        \includegraphics[width=\linewidth]{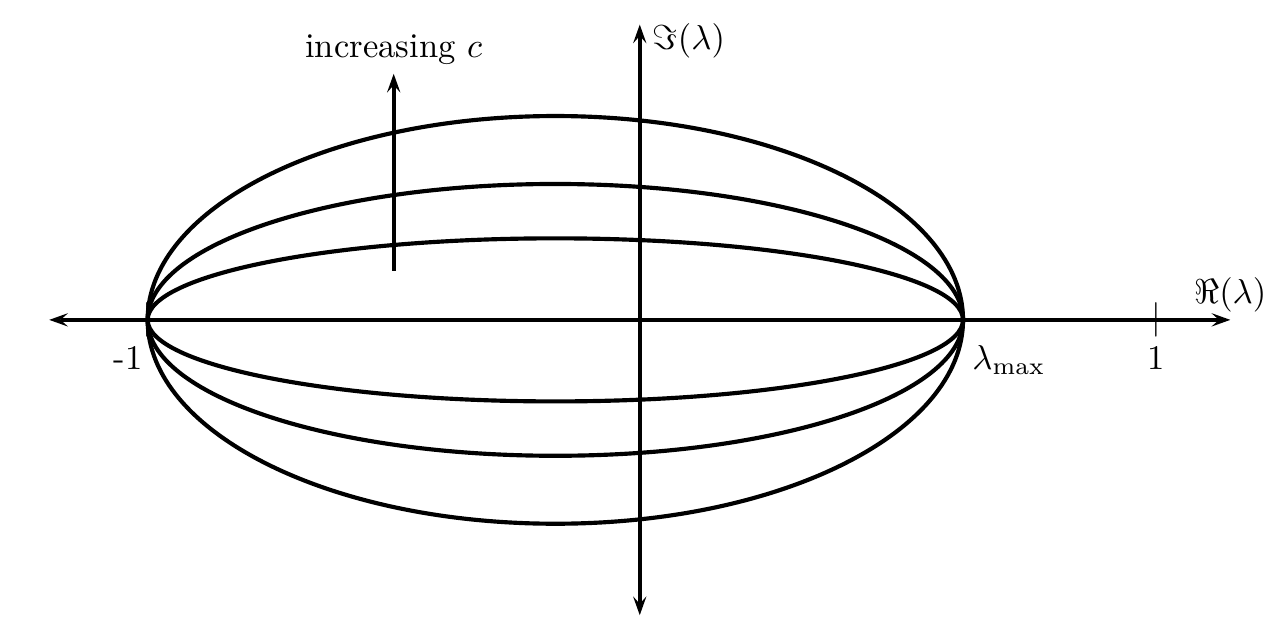}  
        \caption{Varying $c$}
        \label{fig:vary c}
    \end{subfigure}
     \caption{Visualization of the elliptical region over which the minimization problem is solved. For larger $M$, the semi-major axis (range over the real line) increases. For larger $c$, the problem is solved on a thicker ellipse with larger semi-minor axis length.}
    \label{fig:ellipse regions}
\end{figure}

It is clear that the scheme which is the solution to the optimization problem in \eqref{eqn:minimization-srj-real} will only bound the error amplification for $\lambda \in (-1,\lambda_{\text{max}})$. There are no guarantees on the amplification behavior outside this region. To ensure that the schemes developed by solving the optimization problem bound the amplification in the complex plane, we can simply extend the region over which the maximum amplification is minimized. We are particularly interested in developing schemes which attenuate the amplification over elliptical regions of the complex plane, as our non-elliptic PDEs of interest have eigenvalues distributed in this manner. This extension results in the slightly modified optimization problem in Equation \eqref{eqn:optimization-ellipse}.
\begin{equation}
    \min_{\omega_i} \max_{\lambda \in \mathcal{E}} G_{M} (\lambda, \omega_i) \;,\quad\mbox{where}\quad G_M(\lambda;\omega_i):= \prod_{i=1}^{M} \left[(1-\omega_i) + \omega_i \lambda\right]
    \label{eqn:optimization-ellipse}
\end{equation}
In Equation \eqref{eqn:optimization-ellipse}, $\mathcal{E}$ represents the elliptical region of the complex plane over which we wish to minimize the maximum amplification. We specify this elliptical region based on $M$. The ellipse is set to have a semi-major axis which is the same as the portion of the real axis which we previously optimized the scheme over (i.e. $\lambda \in (-1, \lambda_{\text{max}})$) which is specific to each value of $M$. As $M$ increases, the semi-major axis length of the ellipse increases (as depicted in Figure \eqref{fig:vary M}). We consider ellipses corresponding to a few possible semi-minor axis lengths, described by the parameter $c$ which represents the ratio of the semi-minor to semi-major axis. A larger $c$ corresponds to a thicker elliptical region (as depicted in Figure \eqref{fig:vary c}). From these specifications, the elliptical region is described by the following equation \eqref{eqn:ellipse-region} (where $x = \text{Re}(\lambda)$ and $y = \text{Im}(\lambda)$)
\begin{equation}
    \frac{(x-x_{c})^2}{a^2} + \frac{(y)^2}{b^2} \leq 1
    \label{eqn:ellipse-region}
\end{equation}
where 
\begin{align}
    \label{eqn:semi-major}
    a &= \frac{\lambda_{\max} + 1}{2} \\
    \label{eqn:semi-minor}
    b &= c \left( \frac{\lambda_{\max} + 1}{2} \right) \\
    \label{eqn:offset}
    x_{c} &= \frac{\lambda_{\max} - 1}{2} 
\end{align}
In Equation \eqref{eqn:ellipse-region}, $a$ represents the semi-major axis length, $b$ represents the semi-minor axis length, and $c$ represents the ratio of semi-minor to semi-major axis length. Furthermore, $x_{c}$ represents the offset of the center of the semi-major axis from the $y$-axis.  

The optimization problem in Equation \eqref{eqn:optimization-ellipse} can be solved for a variety of different $M$ (number of relaxation parameters) and $c$ (ratio of semi-minor to semi-major axis). The elliptical regions corresponding to these cases is shown pictorially in Figure \ref{fig:ellipse regions}. 

By solving the optimization problem in these different cases, we can establish a family of schemes which utilize different number of relaxation factors $M$ and minimize the amplification over different elliptical regions described by $c$. The different schemes in this family may be more appropriate for solving different linear systems. If the linear system is very stiff (when the Jacobi iteration matrix eigenvalues are very close to 1), it may be more appropriate to utilize a scheme associated with large $M$. If the linear system has eigenvalue spectra which are more spread out in the complex plane, it may be more appropriate to use schemes corresponding to larger $c$. Furthermore, if the eigenvalue spectra does not contain any complex eigenvalues, the schemes corresponding to $c = 0$ (where the elliptical region collapses to the real axis) may be appropriate. In this case, solving the minimization problem should yield the original SRJ schemes described in Section \ref{sec:background-on-srj-schemes}. Knowing the eigenvalue distribution of the linear system can be helpful in determining which scheme is optimal to use, although this information may not always be available, especially for very large systems of equations. 

The minimax optimization problem given in Equation \eqref{eqn:optimization-ellipse} is nontrivial to solve. We present a numerical approach for solving this optimization problem. 

\subsection{Numerical solution to the optimization problem}

We wish to solve the minimax optimization problem given in Equation \eqref{eqn:optimization-ellipse} numerically in order to obtain relaxation schemes which are suitable for solving nonsymmetric linear systems. To make the optimization problem more tractable, we convert the minimax problem into a constrained minimization problem by introducing an auxiliary variable $\bar{g}$ which represents the maximum amplification taken on by the polynomial $G_{M}$. This parameter also represents the bounding value for our amplification polynomial (which was $\frac{1}{3}$ in the case of the original SRJ schemes for symmetric linear systems). Introducing the parameter $\bar{g}$ results in the following constrained minimization problem given in Equation \eqref{eqn:optimization-ellipse-smooth} whose solution is the same as that of the original optimization problem \cite{waren1967}. The optimization problem in Equation \eqref{eqn:optimization-ellipse-smooth} seeks a set of relaxation parameters $\omega_i$ that minimizes the amplification at a collection of test points $x_j$ within our minimization region (these are the optimization constraints). 
\begin{equation}
    \min_{(\omega_i, \bar{g})} \bar{g}^2 \ \ \text{s.t.} \ \ \left| G_{M}(x_j; \omega_i) \right|^2 \leq \bar{g}^2 \ \forall \ x_j \;,\quad\mbox{where}\quad G_M(x_j;\omega_i):= \prod_{i=1}^{M} \left[(1-\omega_i) + \omega_i x_j\right]
    \label{eqn:optimization-ellipse-smooth}
\end{equation}
The test points $x_{j}$ can be chosen arbitrarily, but should ideally be chosen to coincide with the potential extrema locations of the desired amplification polynomial. This ensures that the optimization routine converges towards a solution which bounds the maximum value taken on by the polynomial. For a given $M$, in the case where the optimization is performed over the real line, the locations of the extrema are known. Recall that the amplification polynomials are simply the Chebyshev polynomials under an affine transformation. Therefore, given the extrema locations of the Chebyshev polynomial (denoted by $x^{T_{M}}_{i}$) which are known analytically \cite{mason2002chebyshev}
\begin{equation}
    x^{T_{M}}_{i} = \cos{\left( \frac{i}{M} \pi \right)}, \ \ i = 0, 1, ..., M
\end{equation}
the extrema locations of the amplification polynomials (denoted by $x^{G_{M}}_{i}$) are those of the Chebyshev polynomial under the same transformation (Equation \eqref{eqn:transformation}) as follows
\begin{equation}
    x^{G_{M}}_{i} = g(x^{T_{M}}_{i}) , \ \ i = 0, 1, ..., M \;,\quad\mbox{where}\quad g(\lambda) = \frac{2 \lambda}{1+\lambda^{*}} + \frac{1 - \lambda^{*}}{1+\lambda^{*}}
    \label{eqn:real-test-points}
\end{equation}
where $\lambda^{*}$ satisfies $T_{M}(\lambda^{*}) = 3$. In total, the $M$ order amplification polynomial has $M + 1$ total extrema locations, so we would specify $M+1$ constraints in the optimization problem (one at each extremum location).

\begin{figure}[htbp!]
    \begin{subfigure}{0.5\textwidth}
        \centering
        \includegraphics[width=\linewidth]{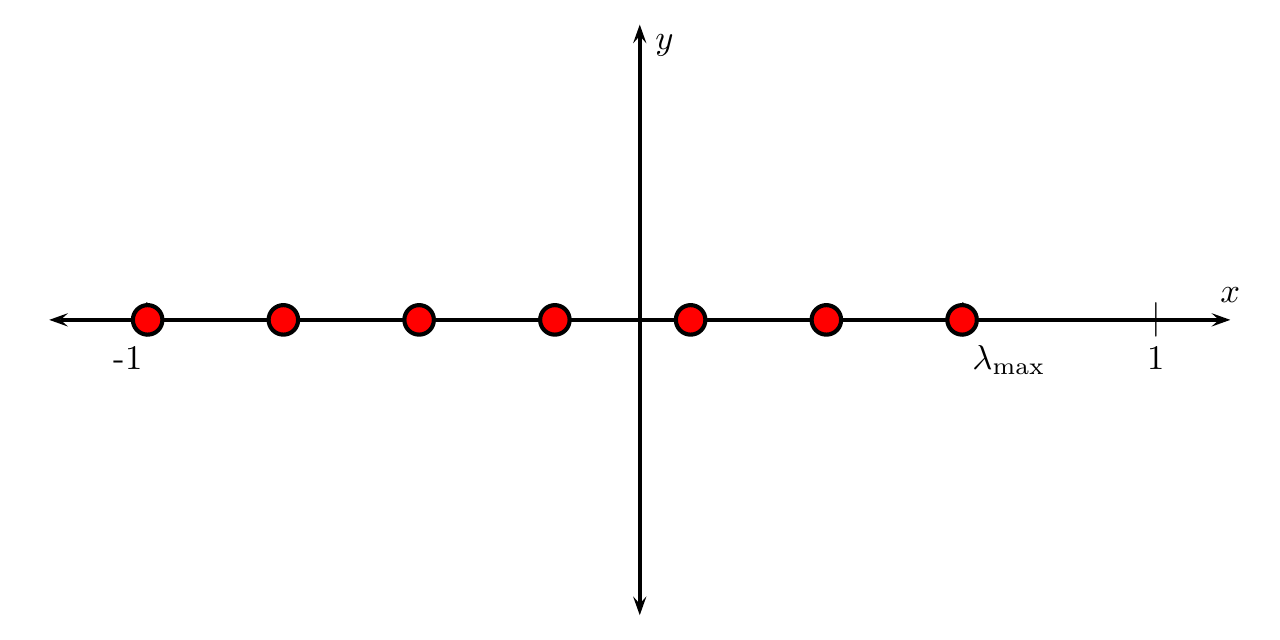}  
        \caption{Real test points}
        \label{fig:real_test_points}
    \end{subfigure}
    \begin{subfigure}{0.5\textwidth}
        \centering
        \includegraphics[width=\linewidth]{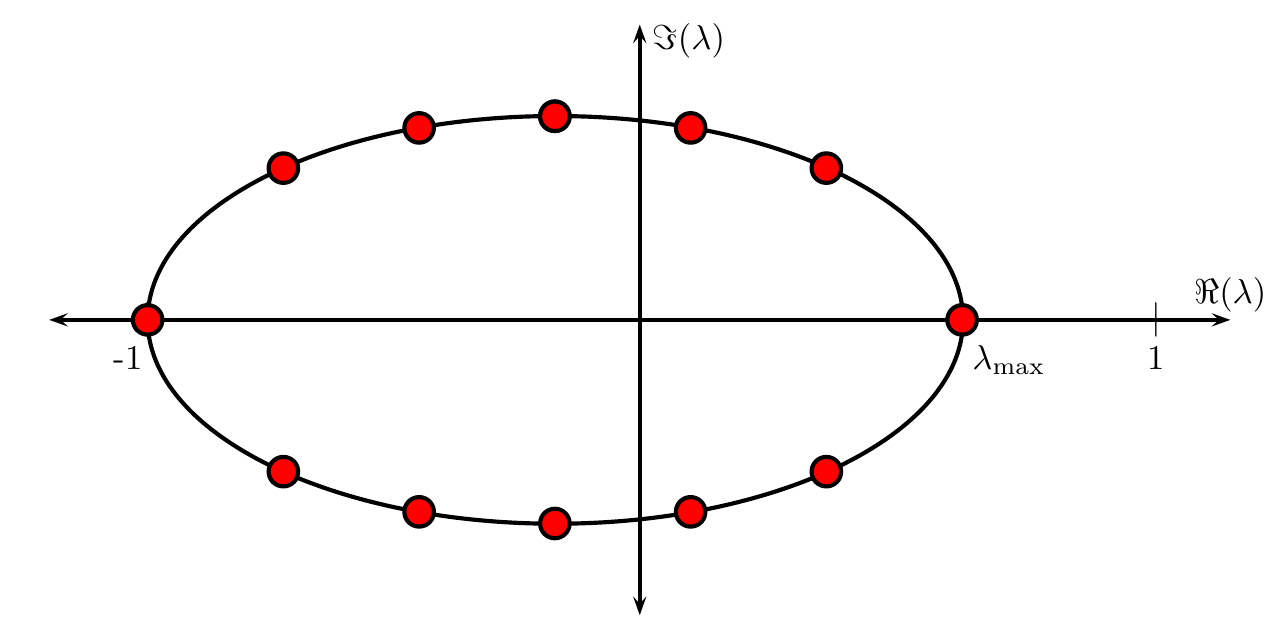}  
        \caption{Complex test points}
        \label{fig:ellipse_test_points}
    \end{subfigure}
    \caption{A visual depiction of the test points over which the constraints are enforced in the optimization problem in Equation \eqref{eqn:optimization-ellipse-smooth}. When minimizing over the real axis for the $M$ order SRJ scheme, $M+1$ constraints are enforced. When minimizing over the elliptical regions, $2M$ constraints are enforced on the ellipse boundary.}
    \label{fig:test_points}
\end{figure}

When performing the optimization over an elliptical region in the complex plane, the test points should be specified on the boundary of the ellipse, as the maximum absolute value of the amplification will be achieved on this boundary (due to the maximum modulus principle from complex analysis \cite{ahlfors1979complex}). In our procedure, we select test points with real value which are the same as those selected in the real case, and with imaginary values such that the test points lie on the ellipse. In other words, each real test point can be mapped to two complex test points for the elliptical case (one on the upper and one on the lower surface of the ellipse), except for the leftmost and rightmost real test points which become mapped to a single real valued test point (see Figure \ref{fig:ellipse_test_points}). Mathematically, these test point locations are expressed as
\begin{equation}
    z_{\text{Test Point}} = x_{j}^{G_{M}} \pm i \left[ b \sqrt{1-\frac{(x_{j}^{G_{M}} - x_{c})^2}{a^2}} \right] , \ \ j = 0, 1, ..., M 
    \label{eqn:ellipse-test-points}
\end{equation}
where $a$, $b$, and $x_{c}$ are specified in Equations \eqref{eqn:semi-major}, \eqref{eqn:semi-minor} and \eqref{eqn:offset}. In total, we can specify the constraints over $2M$ points on the ellipse. A visual depiction of the test points for the real and complex cases is shown in Figure \eqref{fig:test_points}. 

We provide some details regarding the numerical procedure for solving the optimization problem given in Equation \eqref{eqn:optimization-ellipse-smooth}. The $\texttt{scipy.optimize.minimize}$ numerical optimization toolbox is used to solve the constrained optimization problem, using the Trust-Region Constrained numerical algorithm \cite{conn2000trust}. To derive an SRJ scheme of length $M$, our optimization problem has a solution vector of length $M+1$, consisting of the $M$ relaxation factors $w_i$ and amplification bound $\bar{g}$. We denote this solution vector as $\vec{x} = (w_i, \bar{g})$. To solve the optimization problem, we must supply an objective function as well as the Jacobian. The scalar objective function for our problem is given simply by
\begin{equation}
    J(\vec{x}) = \bar{g}^2
    \label{eqn:objective}
\end{equation}
The associated Jacobian vector is the following
\begin{equation}
    \frac{\partial J}{\partial \vec{x}} = [0, \hdots, 0, 2 \bar{g}]^{T}
    \label{eqn:objective-jacobian}
\end{equation}
and is simply a vector of length $M+1$ with zero as the first $M$ entries, and $2 \bar{g}$ in the last entry. We must also specify the constraint functions for our constrained optimization problem, as well as any Jacobians/Hessians associated with these constraints. We write the constraint in the form of an inequality $c(\vec{x}) \geq 0$ where $c(\vec{x})$ is given by
\begin{equation}
    c(\vec{x}) = \bar{g}^2 - \left| G_{M}(x_j; \omega_i) \right|^{2} 
    \label{eqn:constraint}
\end{equation}
The corresponding Jacobian vector for the constraint function has the form
\begin{equation}
    J = \frac{\partial c(\vec{x})}{\partial \vec{x}} = \bigg[\underbrace{\frac{\partial c(\vec{x})}{\partial \omega_j}}_{M \ \text{terms}}, \underbrace{\frac{\partial c(\vec{x})}{\partial \bar{g}}}_{1 \ \text{term}}\bigg]^{T}
    \label{eqn:jacobian-vector-constraint}
\end{equation}
and is comprised of $M+1$ terms corresponding to the derivative of the constraint function with respect to the $M$ relaxation factors and the bounding value $\bar{g}$ which comprise the solution variables in $\vec{x}$.
A derivation of the Jacobian vector entries and a procedure for constructing it is provided in Appendix \ref{sec:appendix1-jacobian-derivation}. Lastly, a 3-point finite difference approximation is used to compute the constraint Hessians based on the constraint Jacobian.

We solve the optimization problem in Equation \eqref{eqn:optimization-ellipse-smooth} with the associated objective function and objective Jacobian vector in Equations \eqref{eqn:objective} and \eqref{eqn:objective-jacobian}, and associated constraints and constraint Jacobian in Equations \eqref{eqn:constraint} and \eqref{eqn:jacobian-vector-constraint}. The constraint, constraint Jacobian and constraint Hessian are specified at all test points (given by Equation \eqref{eqn:real-test-points} for the real case or Equation \eqref{eqn:ellipse-test-points} for the ellipse case) for the optimization routine. SRJ schemes which are solutions to the optimization problem are found for $M = 2$ up to $M = 20$ relaxation parameters. For each $M$, several schemes are obtained resulting from an optimization over ellipses of varying thickness corresponding to $c = 0, 1/10, 1/5, 1/3, 1/2$. The schemes resulting from $M = 2$ for the various elliptical regions corresponding to different $c$ are shown in Table \ref{tab:srj-schemes-M2}. Schemes corresponding to $M = 5$ are also shown in Table \ref{tab:srj-schemes-M5}. The schemes corresponding to $c = 0$ are the original SRJ schemes derived for the symmetric case. The full set of schemes derived in this work is provided in Appendix \ref{sec:all-srj-schemes}.

\begin{figure}[htbp!]
    \centering
    \includegraphics[width=0.7\textwidth]{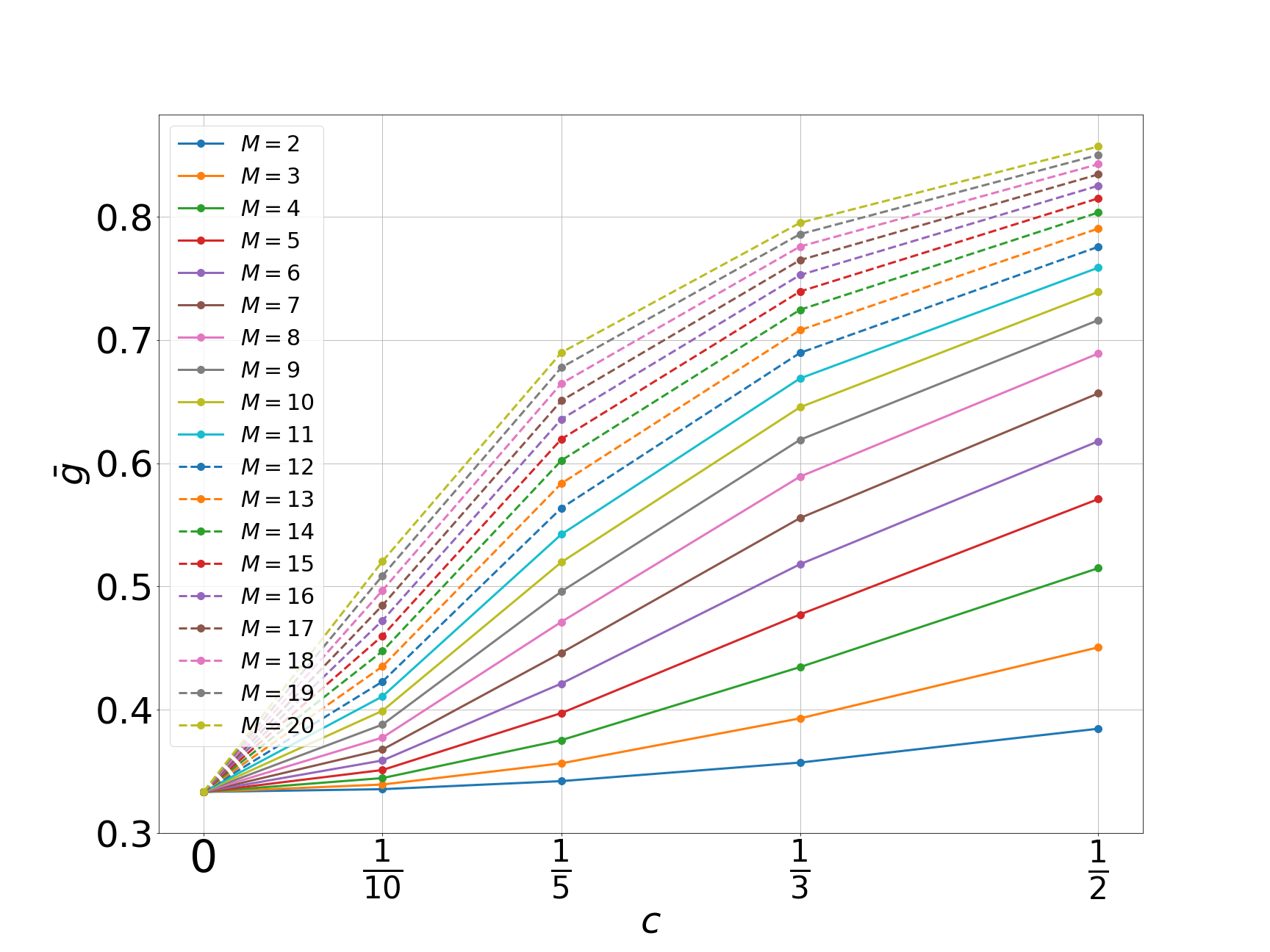}
    \caption{Maximum amplification bounding value $\bar{g}$ corresponding to each SRJ scheme}
    \label{fig:optimization-robustness}
\end{figure}

The maximum amplification or bounding value corresponding to each SRJ scheme for each combination of $M$ and $c$ (obtained from the optimization routine as $\bar{g}$) is summarized in Figure \ref{fig:optimization-robustness}. As $M$ and $c$ increase, the resulting maximum amplification has a greater magnitude. This is expected, as increasing $M$ or $c$ effectively increases the size of the elliptical region over which the minimization problem is solved. Figure \ref{fig:optimization-robustness} shows that the maximum amplification associated with each scheme has magnitude less than 1 within the elliptical region corresponding to its $M$ and $c$ value. This indicates that the schemes are expected to converge if the iteration matrix eigenvalues of the linear system of interest are contained within the elliptical region over which the scheme minimizes the amplification. 

\begin{table}[htbp!]
    \caption{Relaxation factors associated with SRJ schemes for $M = 2$ and various ellipse semi-minor to semi-major length ratios}
    \centering
    \begin{tabular}{|c|l|}
        \hline
        $c$ & SRJ scheme parameters \\
        \hline
        0 & 1.70710678, 0.56903559 \\
        \hline
        $1/10$ & 1.6985778 , 0.56998629 \\
        \hline
        $1/5$ & 1.67329915, 0.57289385\\
        \hline
        $1/3$ & 1.61475726, 0.58009431 \\
        \hline
        $1/2$ & 1.50541883, 0.59563558 \\
        \hline
    \end{tabular}
    \label{tab:srj-schemes-M2}
\end{table}

\begin{table}[htbp!]
    \caption{Relaxation factors associated with SRJ schemes for $M = 5$ and various ellipse semi-minor to semi-major length ratios}
    \centering
    \begin{tabular}{|c|l|}
        \hline
        $c$ & SRJ scheme parameters \\
        \hline
        0 & 9.23070078, 0.62486986, 0.51215172, 2.17132939, 0.97045898\\
        \hline
        $1/10$ & 2.15794366, 8.85298329, 0.62598725, 0.51336697, 0.9704587\\
        \hline
        $1/5$ & 0.97045888, 2.11836786, 7.87621951, 0.62939827, 0.51708554\\
        \hline
        $1/3$ & 2.02782132, 0.52636836, 0.97045899, 6.20847021, 0.63786058 \\
        \hline
        $1/2$ & 0.54674459, 4.31270705, 0.65617569, 0.97045902, 1.86254896\\
        \hline
    \end{tabular}
    \label{tab:srj-schemes-M5}
\end{table}

\begin{figure}[htbp!]
    \begin{subfigure}{.5\textwidth}
        \centering
        \includegraphics[width=\textwidth]{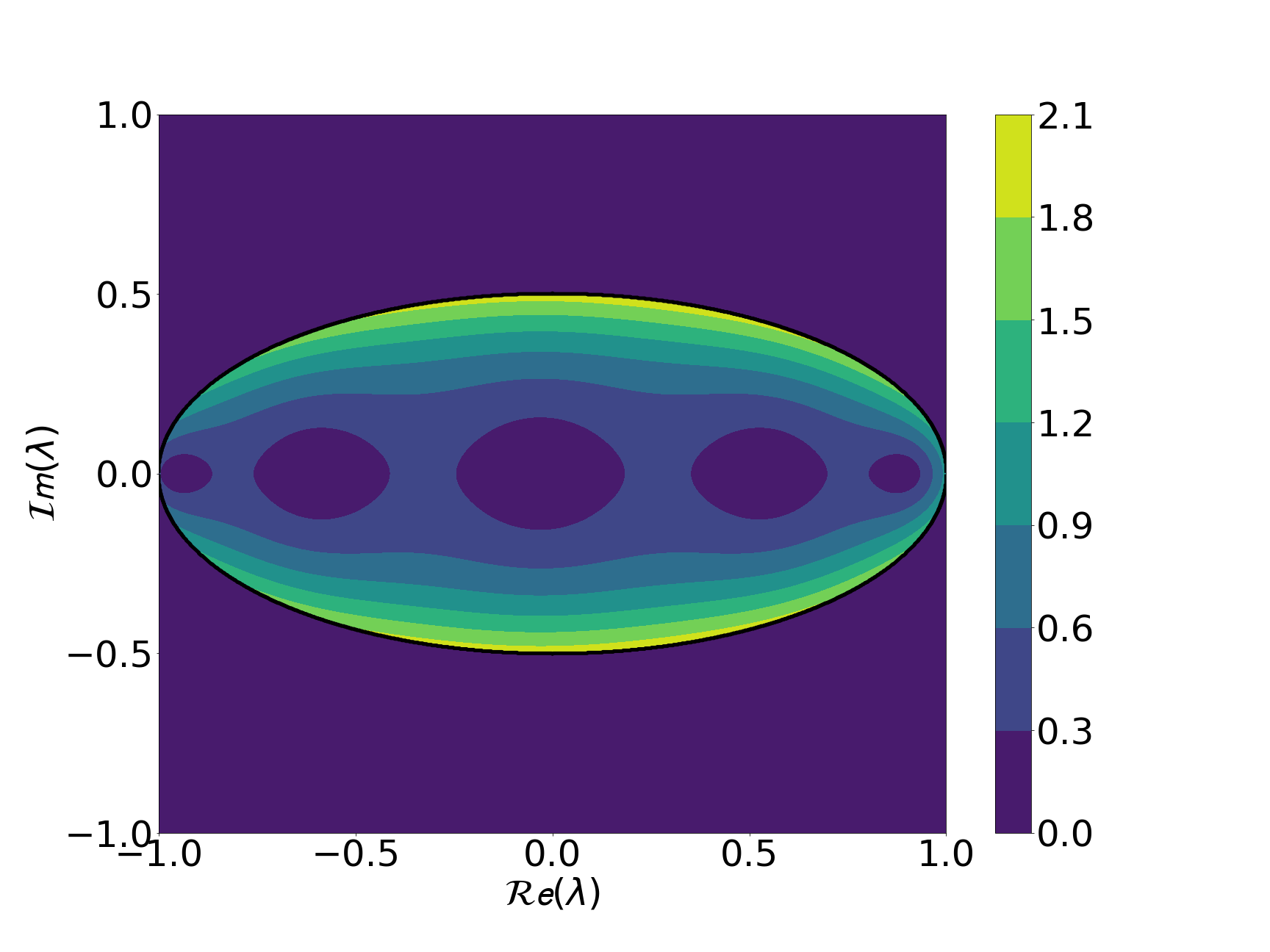}
        \caption{Original $M = 5$ SRJ amplification}
        \label{fig:ellipse_amplification_M5_original}
    \end{subfigure}
    \begin{subfigure}{.5\textwidth}
        \centering
        \includegraphics[width=\textwidth]{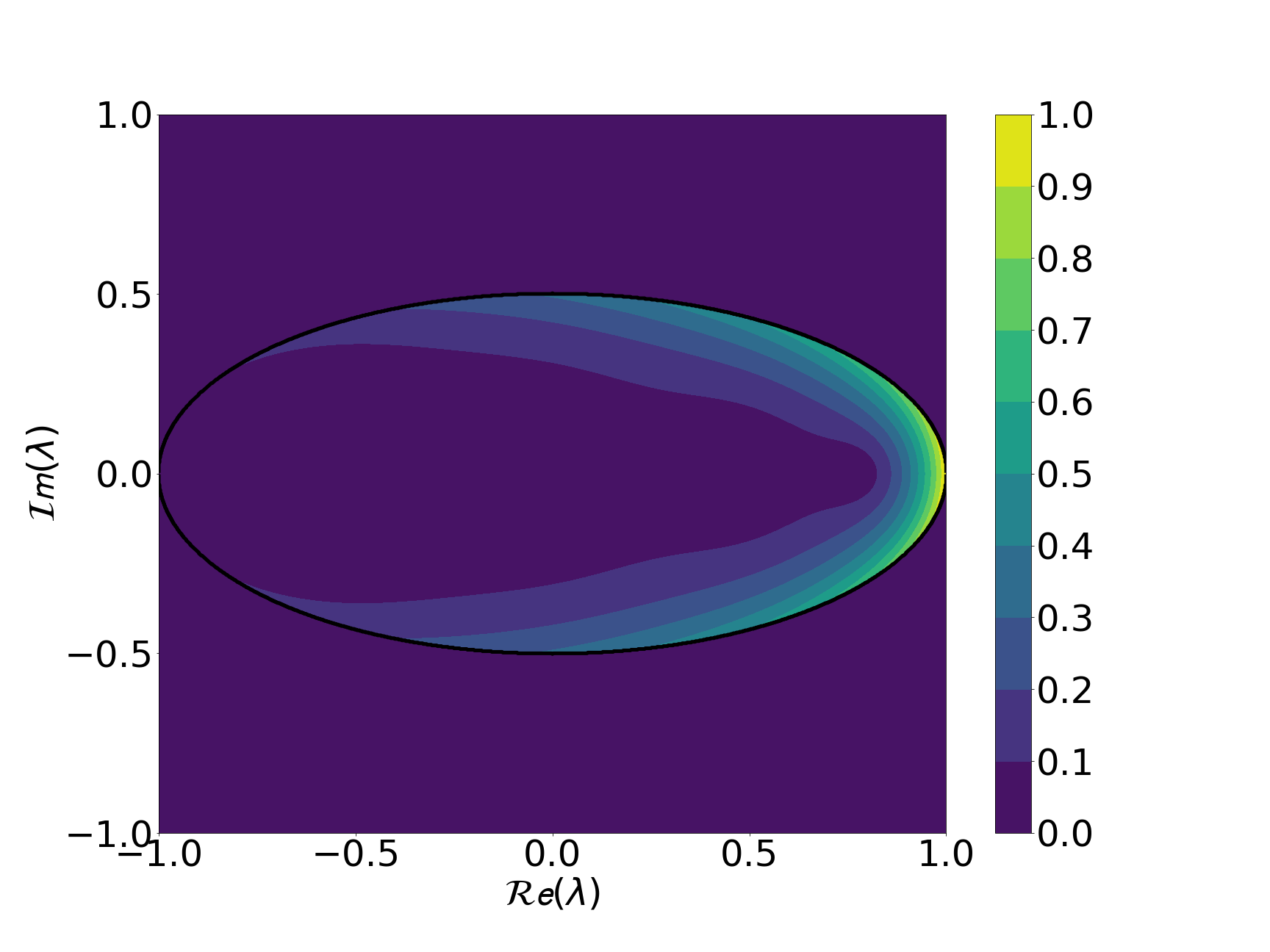}
        \caption{Optimized $M = 5$ SRJ amplification}
        \label{fig:ellipse_amplification_M5_optimized}
    \end{subfigure}
    \\
    \begin{subfigure}{.5\textwidth}
        \centering
        \includegraphics[width=\textwidth]{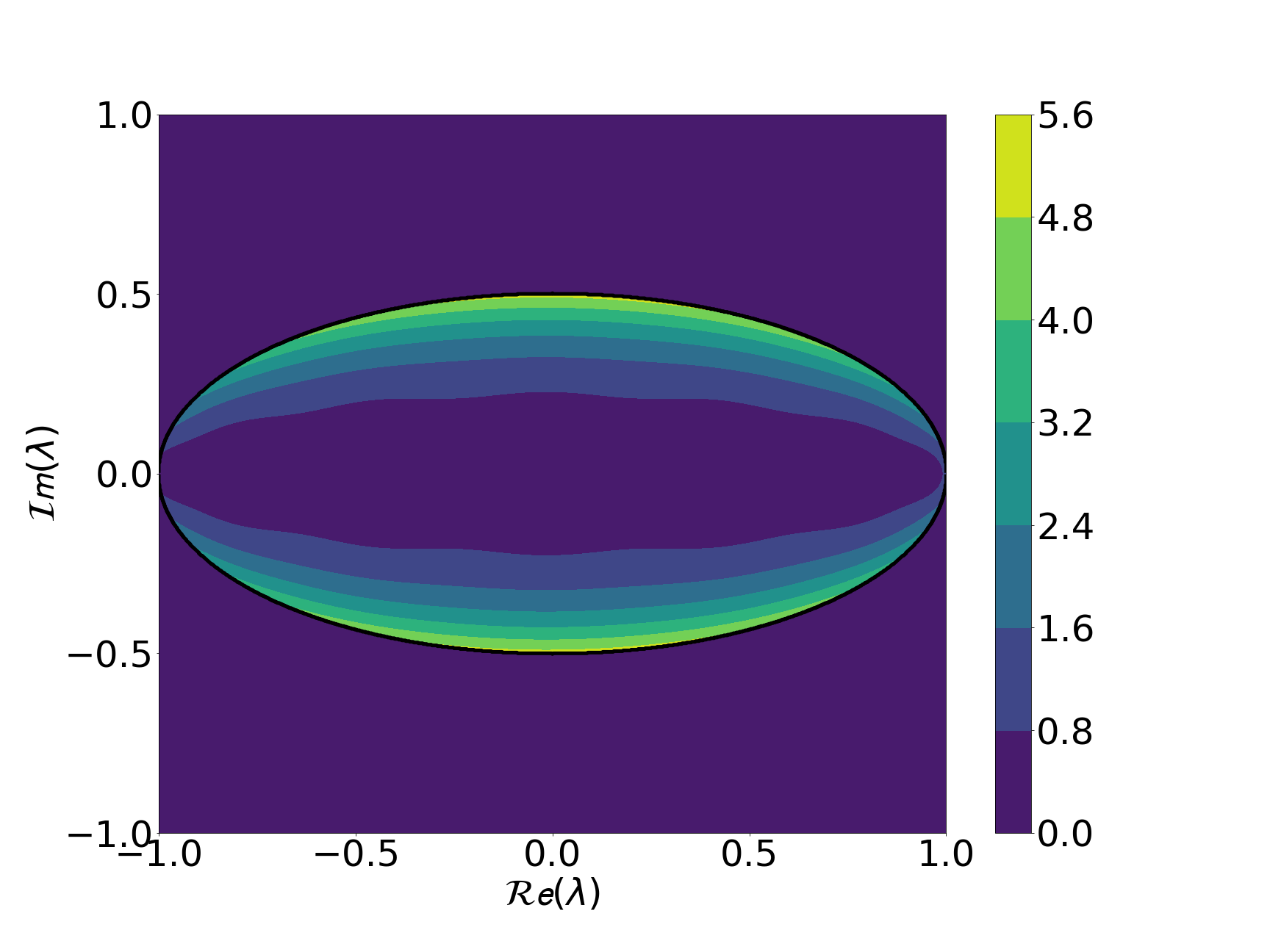}
        \caption{Original $M = 7$ SRJ amplification}
        \label{fig:ellipse_amplification_M5_original}
    \end{subfigure}
    \begin{subfigure}{.5\textwidth}
        \centering
        \includegraphics[width=\textwidth]{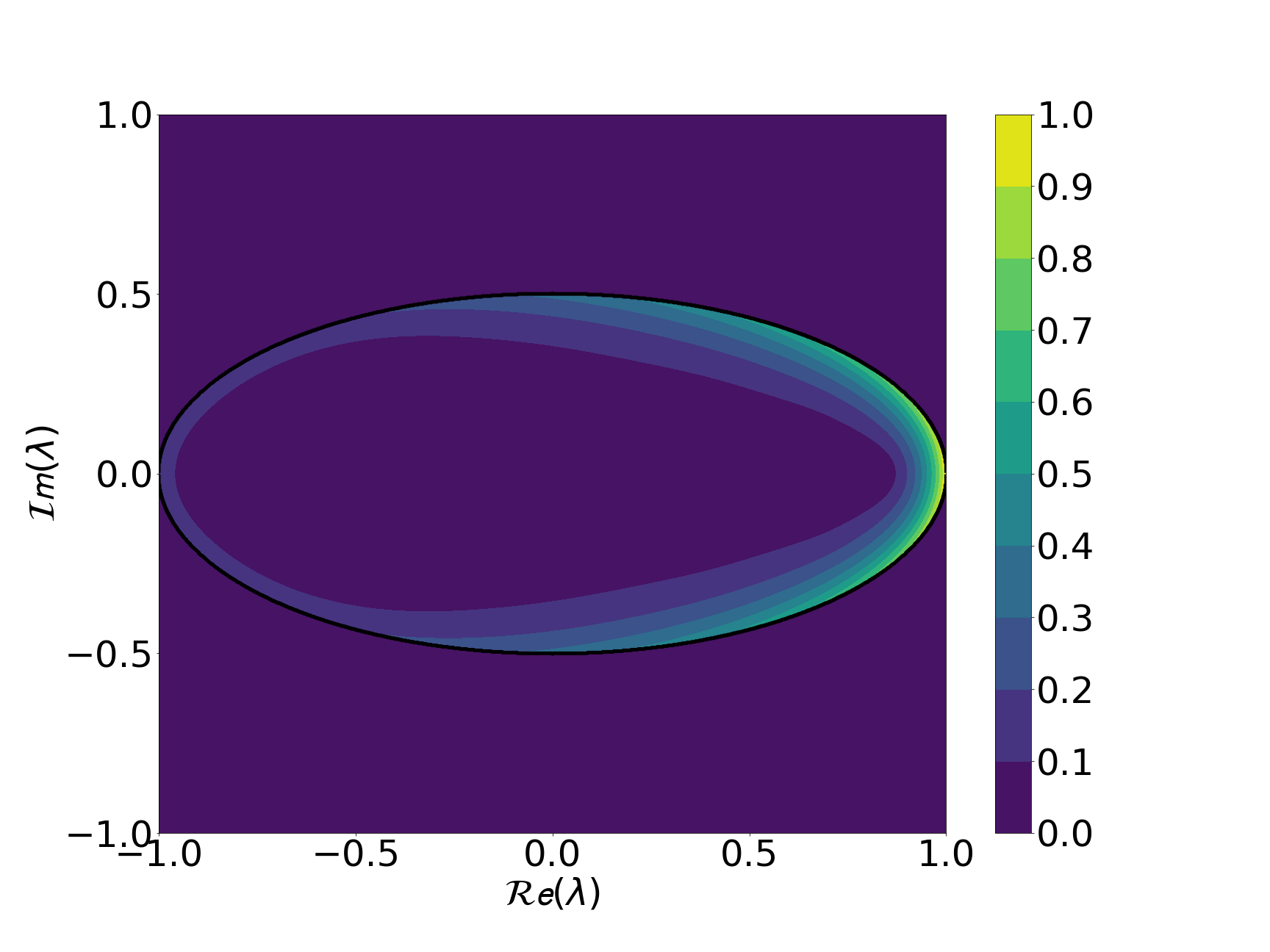}
        \caption{Optimized $M = 7$ SRJ amplification}
        \label{fig:ellipse_amplification_M5_optimized}
    \end{subfigure}
    \caption{Comparison of the amplification due to the original SRJ schemes (left) and newly developed SRJ schemes (right) optimized over ellipse corresponding to $c = \frac{1}{2}$ and specified $M$. The amplification of the original schemes exceeds 1 in this region, but is bounded for the new schemes.}
    \label{fig:amplification_comparison_ellipse}
\end{figure}

To illustrate the benefit of the new schemes, we compare the expected amplification behavior of the original $M = 5, 7$ SRJ schemes to the new $M = 5,7$ schemes optimized over an ellipse region corresponding to $c = 1/2$. The amplification for all schemes is shown in this $c = 1/2$ region corresponding to the given $M$. We observe that the original scheme amplification (shown on the left in Figure \eqref{fig:amplification_comparison_ellipse}) exceeds 1 near the boundary of this region. However, the optimized schemes are bounded by 1 here (shown on the right in Figure \eqref{fig:amplification_comparison_ellipse}). If the schemes were used to solve a linear system whose iteration matrix eigenvalues lie close to the ellipse boundaries, the original schemes are very likely to diverge while the new schemes would converge. For even larger $M$, the maximum amplification achieved in this ellipse region is expected to be larger. Therefore, the original schemes become more unsuitable as many relaxation parameters are used. 

Another metric of comparison for the SRJ schemes derived here are their ability to solve extremely stiff systems. An example of a stiff linear system is the linear system arising from discretization of the 1D Poisson equation on an extremely refined grid. In this case, the spectral radius of the Jacobi iteration matrix is very close to 1, meaning the convergence rate of Jacobi iteration is very slow. One measure of an SRJ scheme's ability to solve stiff systems is the slope of the corresponding amplification polynomial along the real axis at a value of $\lambda = 1$. The spectral radius of the SRJ iteration matrix for a stiff linear system is more likely to deviate from 1 if the slope of the amplification polynomial is greater. The slope associated with each of the schemes derived in this work is listed in Appendix \ref{sec:scheme-slope}. As $M$ increases, the slope is greater indicating that the scheme is more well suited to handling stiffer systems. However, as $c$ increases the slope decreases, indicating that it has more difficulty handling a stiffer system. This presents a tradeoff when selecting a scheme to use to solve a linear system. For stiff systems, a scheme associated with large $M$ and smaller $c$ may be desirable. However, when the iteration matrix eigenvalue spectrum has many complex eigenvalues with a large spread, it may be preferable to use a scheme corresponding to larger $c$. 

In this section, we have presented a methodology for deriving SRJ schemes as the solution to an optimization problem. We have presented evidence that the schemes are expected to work well for solving linear systems with complex eigenvalues (nonsymmetric matrices). The next section illustrates the convergence behavior of these schemes for solving nonsymmetric linear systems arising from non-elliptic PDEs.

\section{Results}

In this section, we demonstrate the convergence behavior of several of the SRJ schemes derived in Section \ref{sec:srj-for-nonsymmetric-schemes}. These schemes are used to solve the one-dimensional and two-dimensional steady advection diffusion equations which are parabolic in nature. The convergence behavior is compared to the standard Jacobi iteration.

\subsection{1D Steady Advection-Diffusion Equation}

We consider the one-dimensional steady advection-diffusion equation given in Equation \eqref{eqn:1d-advection-diffusion}, on a domain $x \in [0,1]$. A homogenous Dirichlet boundary condition is specified on the left of the domain, and a homogenous Neumann boundary condition is specified on the right. 
\begin{equation}
    - \nu \frac{d^{2} u(x)}{dx^2} + a \frac{du(x)}{dx} = f(x), \ u(x = 0) = 0, \ \frac{du}{dx}\bigg|_{x = 1} = 0,
    \label{eqn:1d-advection-diffusion}
\end{equation}
In Equation \eqref{eqn:1d-advection-diffusion}, $\nu$ refers to the diffusion coefficient, $a$ is the advection coefficient, and $f(x)$ is some forcing function. We consider constant advection and diffusion coefficients. The forcing function is set to be $f(x) = \sin(2\pi x)$. To derive a set of linear equations, we discretize the PDE using a second order central difference scheme for the second derivative
\begin{align}
    \frac{d^{2} u}{dx^2} \approx \frac{u_{i+1} - 2u_{i} + u_{i-1}}{\Delta x^2} 
\end{align}
and an upwinding scheme for the first derivative as follows
\begin{equation}
    \frac{du}{dx} \approx
    \begin{cases}
        \frac{u_{i} - u_{i-1}}{\Delta x} \ \text{if} \ a > 0 \\
        \frac{u_{i+1} - u_{i}}{\Delta x} \ \text{if} \ a < 0
    \end{cases}
\end{equation}
where $u_{i}$ is the solution at the $i$th grid point and $\Delta x$ is the grid spacing (we assume a uniform grid). The finite difference discretization above leads to the tridiagonal system of equations $Ax = b$ for the unknowns at the grid points, where $A$ is given below
\[
A = 
\begin{pmatrix}
a_{0} & a_{1} &  \\
a_{-1} & a_{0} & a_{1}  \\
 & \ddots & \ddots & \ddots  \\
 &  & a_{-1} & a_{0} & a_{1}  \\
 &  &  & a_{-1} & a_{0}
\end{pmatrix} 
\] 
with the following entries (assuming $a > 0$)
\begin{align*}
    a_{-1} &= -\frac{\nu}{\Delta x^2} - \frac{a}{\Delta x} \\
    a_{0} &= \frac{2 \nu}{\Delta x^2} + \frac{a}{\Delta x} \\
    a_{1} &= -\frac{\nu}{\Delta x^2}
\end{align*}
except for the last row which must be adjusted to reflect the Neumann boundary condition.
Note that setting $\nu = 1$ and $a = 0$ results in the tridiagonal and symmetric linear system corresponding to the Poisson equation. As the value of $a$ (amount of advection) increases, the matrix becomes increasingly nonsymmetric.

\begin{figure}[htbp!]
    \centering
    \includegraphics[width=0.7\textwidth]{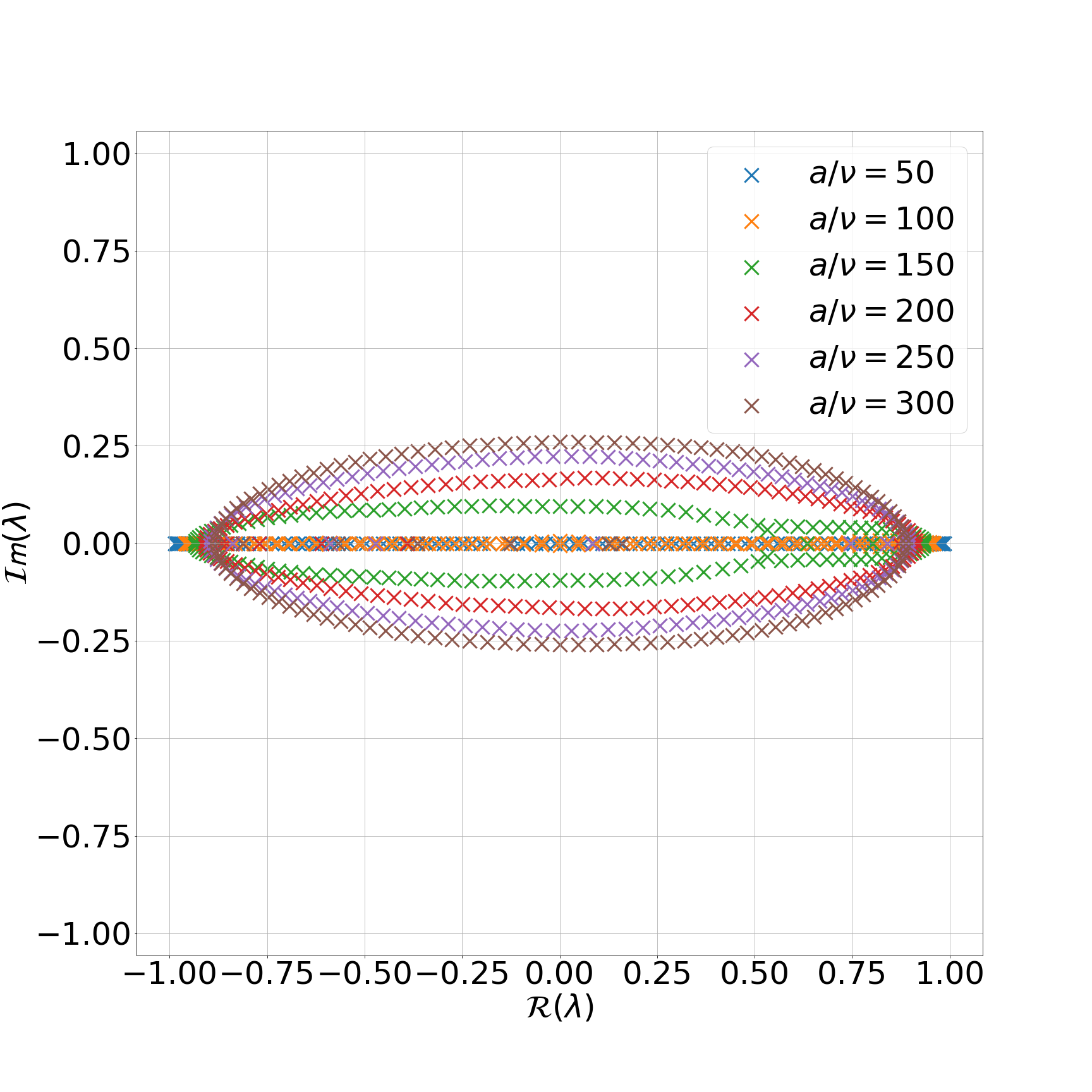}
    \caption{Eigenvalue spectrum of the iteration matrix $B_{\text{J}} = -D^{-1}(L+U)$ for different $\frac{a}{\nu}$ ratios. As $\frac{a}{\nu}$ increases, the eigenvalues encompass a larger region of the complex plane which resembles an ellipse.}
    \label{fig:eigenvalue_spectrum_BJ}
\end{figure}

We test the performance of our SRJ schemes on the linear system arising from a finite difference discretization of this one-dimensional steady advection diffusion equation with $N = 128$ unknowns. For our test problems, we set the diffusion coefficient $\nu = 1$ and vary the advection coefficient according to the following values: $a = 50, 100, 150, 200, 250, 300$. Varying $a$ changes the eigenvalue spectrum of the iteration matrix $B_\text{J}$, which is shown in Figure \ref{fig:eigenvalue_spectrum_BJ} for the specified $a$ values. In particular, as the amount of advection increases, the imaginary components of the eigenvalues of the iteration matrix increase as well. The eigenvalue spectra appear to encompass an elliptical domain, similar to the regions over which the SRJ schemes minimize the amplification in our optimization problem in Equation \eqref{eqn:optimization-ellipse}. These elliptical regions grow thicker as the amount of advection $a$ is increased.

\begin{figure}[htbp!]
    \begin{subfigure}{.5\textwidth}
        \centering
        \includegraphics[width=\textwidth]{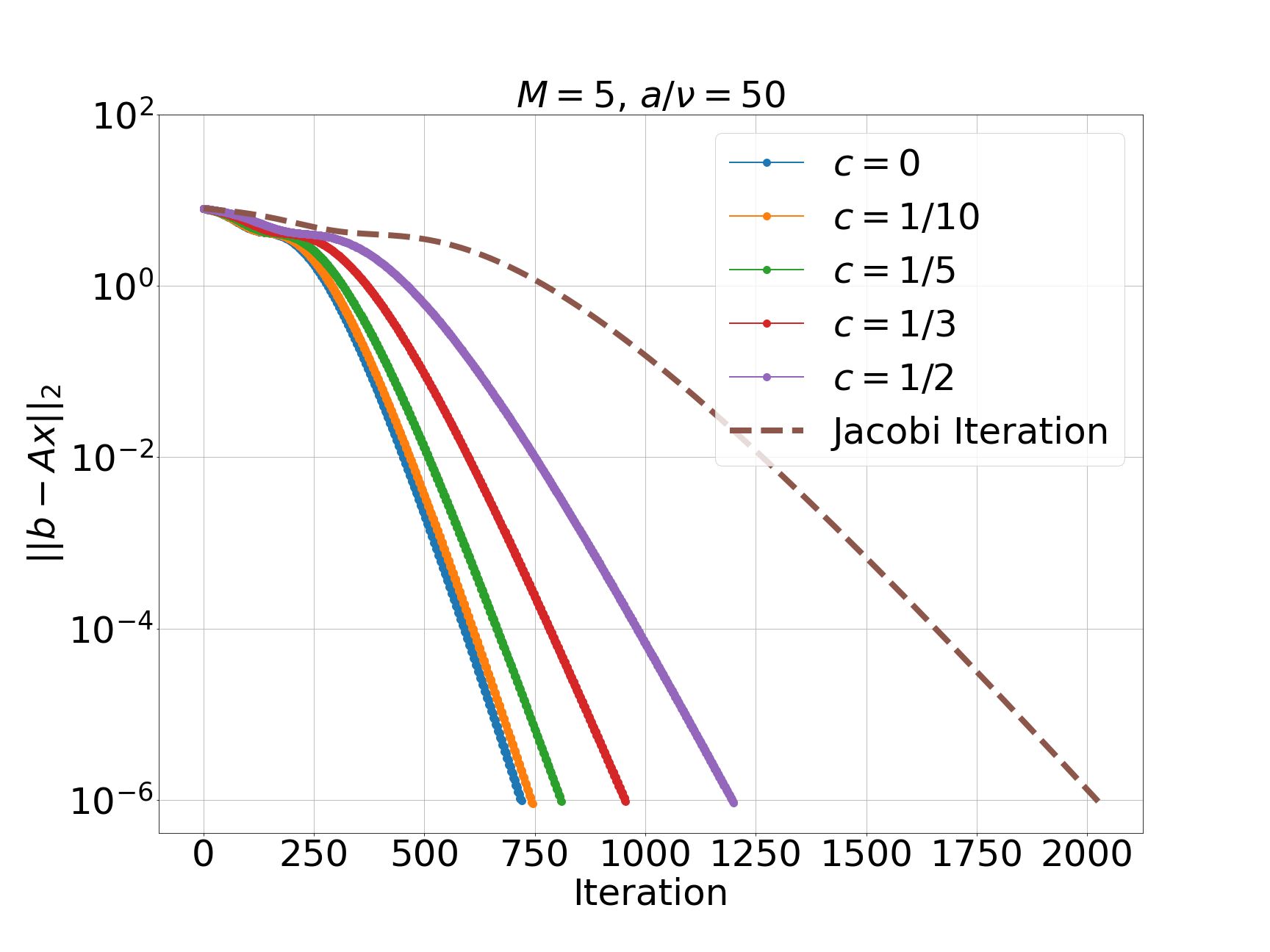}
        \caption{Convergence for $a = 50$}
        \label{fig:convergence_a_50}
    \end{subfigure}
    \begin{subfigure}{.5\textwidth}
        \centering
        \includegraphics[width=\textwidth]{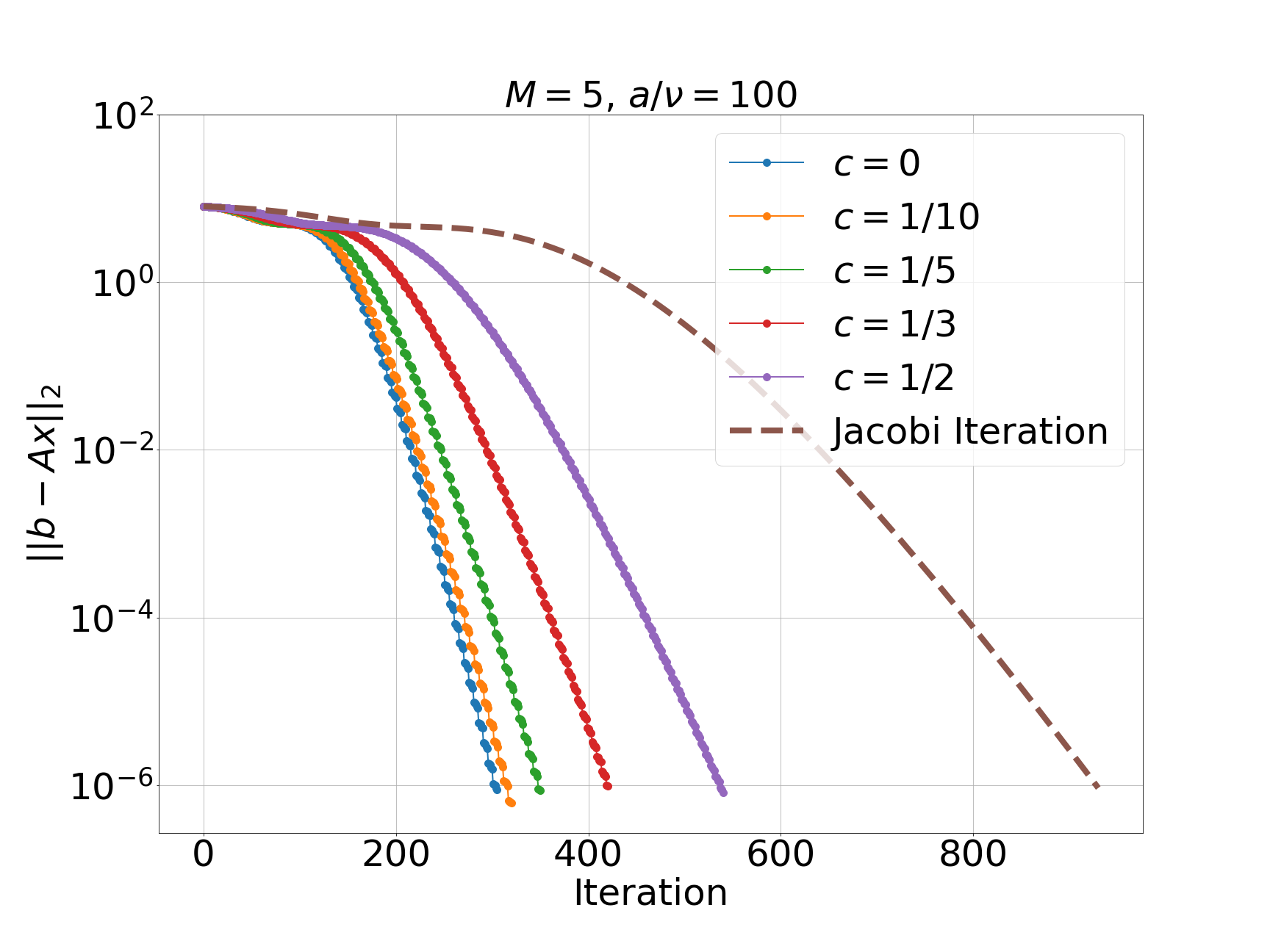}
        \caption{Convergence for $a = 100$}
        \label{fig:convergence_a_100}
    \end{subfigure}
    \\
    \begin{subfigure}{.5\textwidth}
        \centering
        \includegraphics[width=\textwidth]{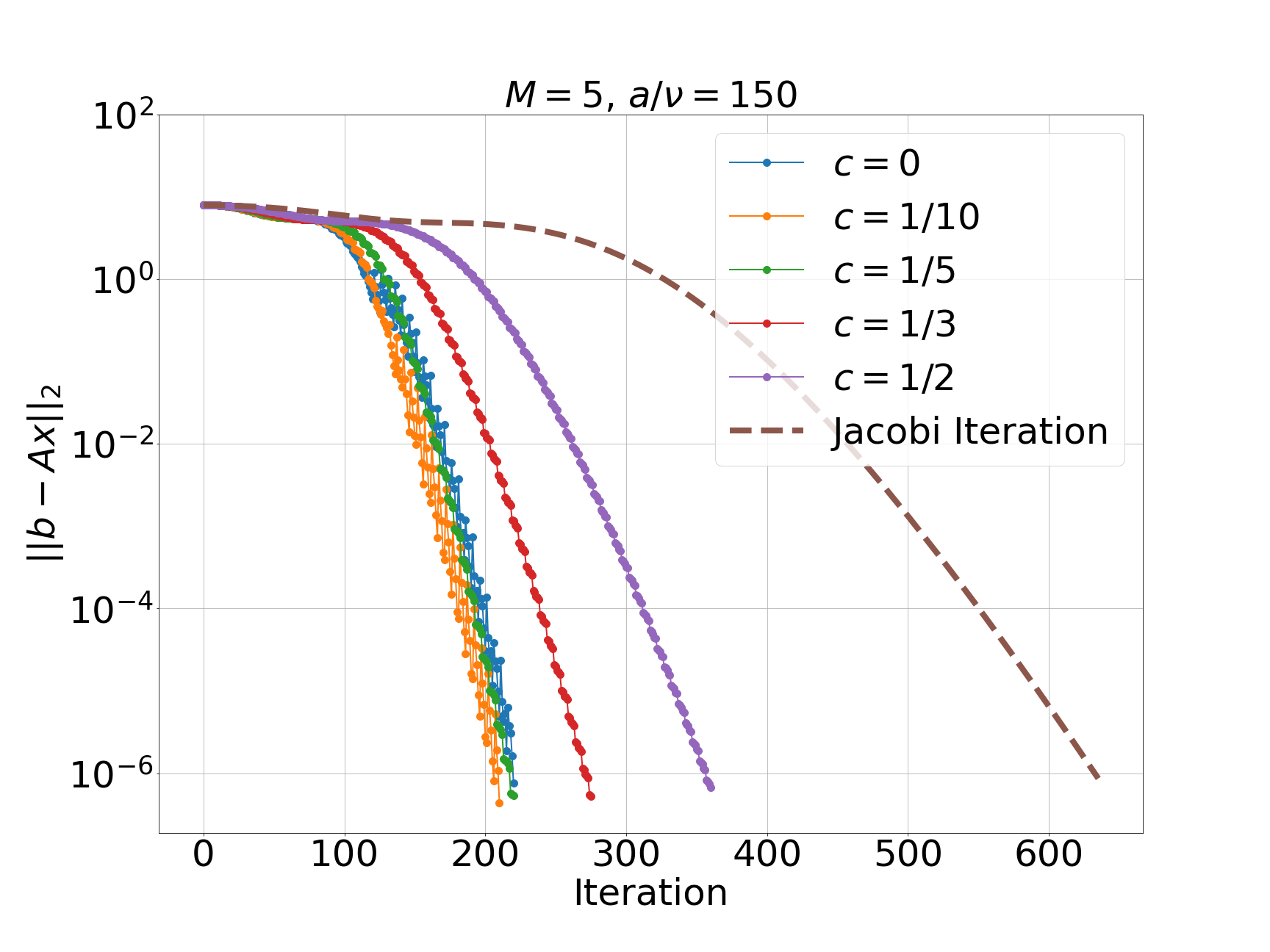}
        \caption{Convergence for $a = 150$}
        \label{fig:convergence_a_150}
    \end{subfigure}
    \begin{subfigure}{.5\textwidth}
        \centering
        \includegraphics[width=\textwidth]{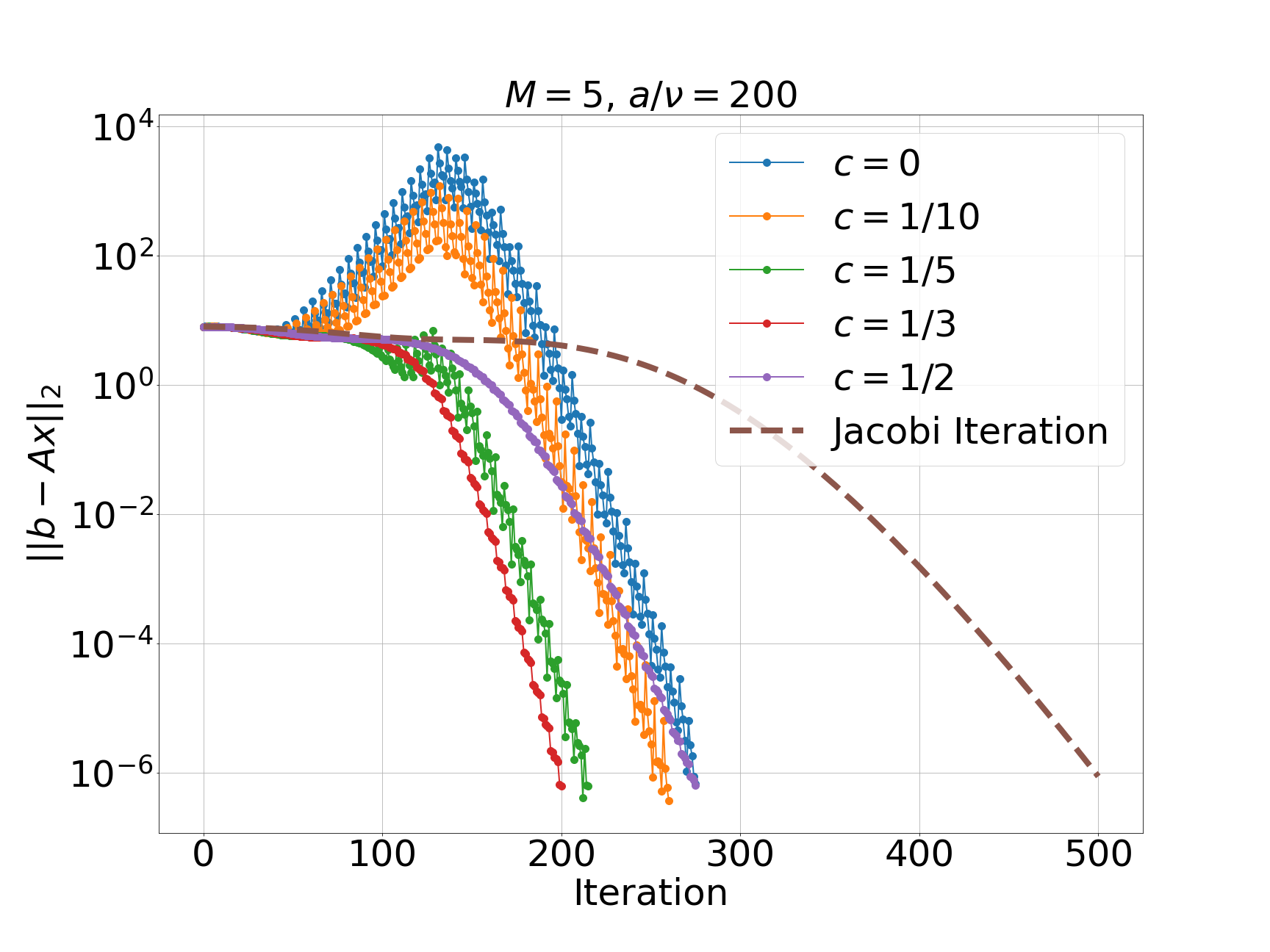}
        \caption{Convergence for $a = 200$}
        \label{fig:convergence_a_200}
    \end{subfigure}
    \\
    \begin{subfigure}{.5\textwidth}
        \centering
        \includegraphics[width=\textwidth]{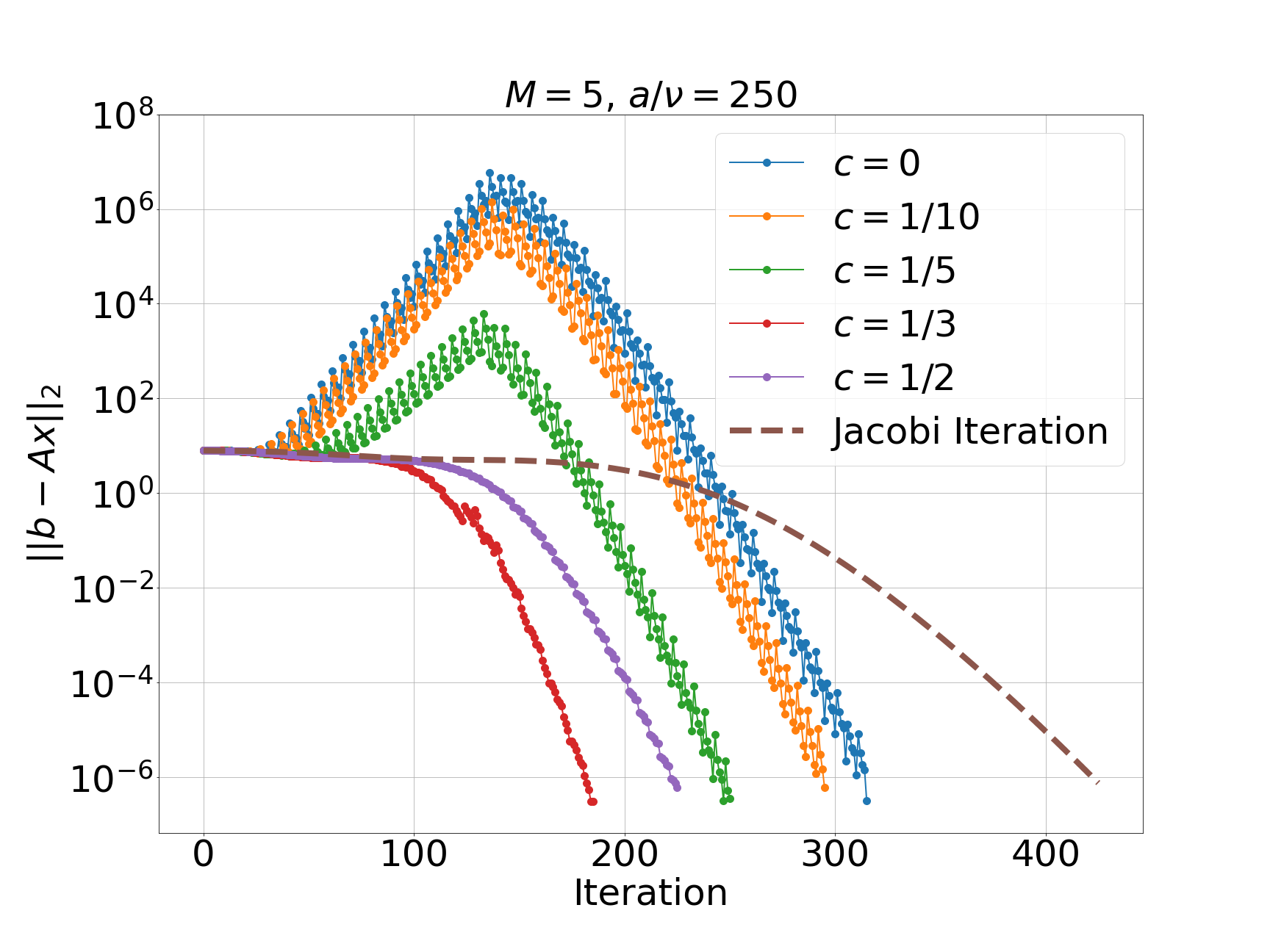}
        \caption{Convergence for $a = 250$}
        \label{fig:convergence_a_250}
    \end{subfigure}
    \begin{subfigure}{.5\textwidth}
            \centering
            \includegraphics[width=\textwidth]{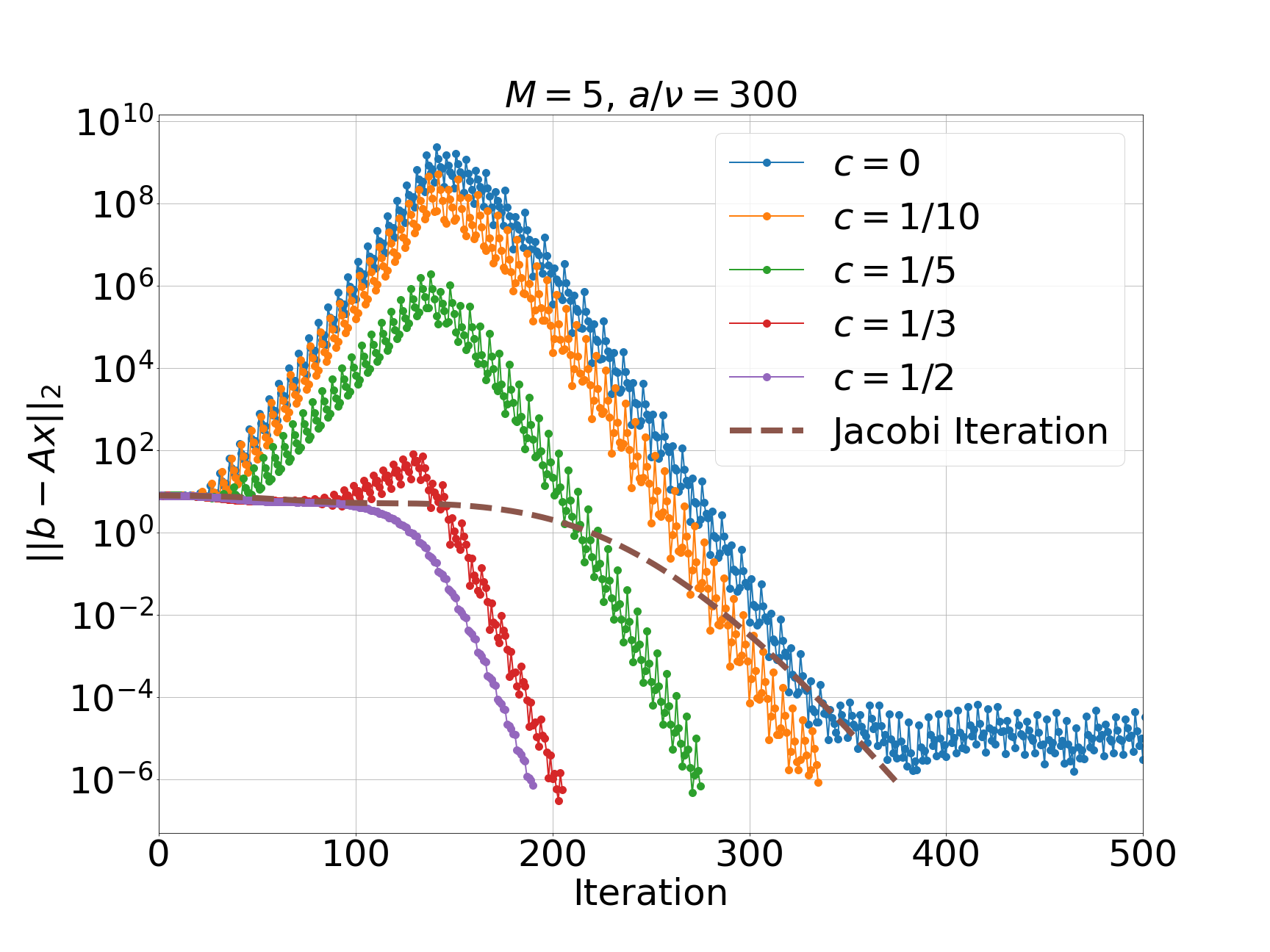}
            \caption{Convergence for $a = 300$}
            \label{fig:convergence_a_300}
    \end{subfigure}
    \caption{Convergence behavior of the $M = 5$ SRJ schemes for solving the 1D advection-diffusion of varying advection to diffusion ratio. As the amount of advection increases, the schemes optimized over a larger elliptical region (associated with larger $c$) exhibit faster convergence relative to the standard Jacobi iteration and schemes optimized only over the real axis ($c = 0$).}
    \label{fig:convergence_of_advection_diffusion}
\end{figure}

We solve six different linear systems, each corresponding to a different advection value $a$. The SRJ schemes corresponding to five relaxation factors ($M = 5$) and varying elliptical thickness $c$ are used to solve each linear system. The convergence of the $M = 5$ schemes associated with different values of the ellipse ratio $c$ are compared and shown in Figure \ref{fig:convergence_of_advection_diffusion}. The convergence associated with the standard Jacobi iteration is also shown. In all cases, the initial solution is set to a vector of ones. The convergence plots show the $L_{2}$ residual norm $||b-Ax||_{2}$ at each SRJ iteration until a residual below a tolerance value of 1E-6 is achieved.

\begin{figure}[htbp!]
    \centering
    \includegraphics[width=0.7\textwidth]{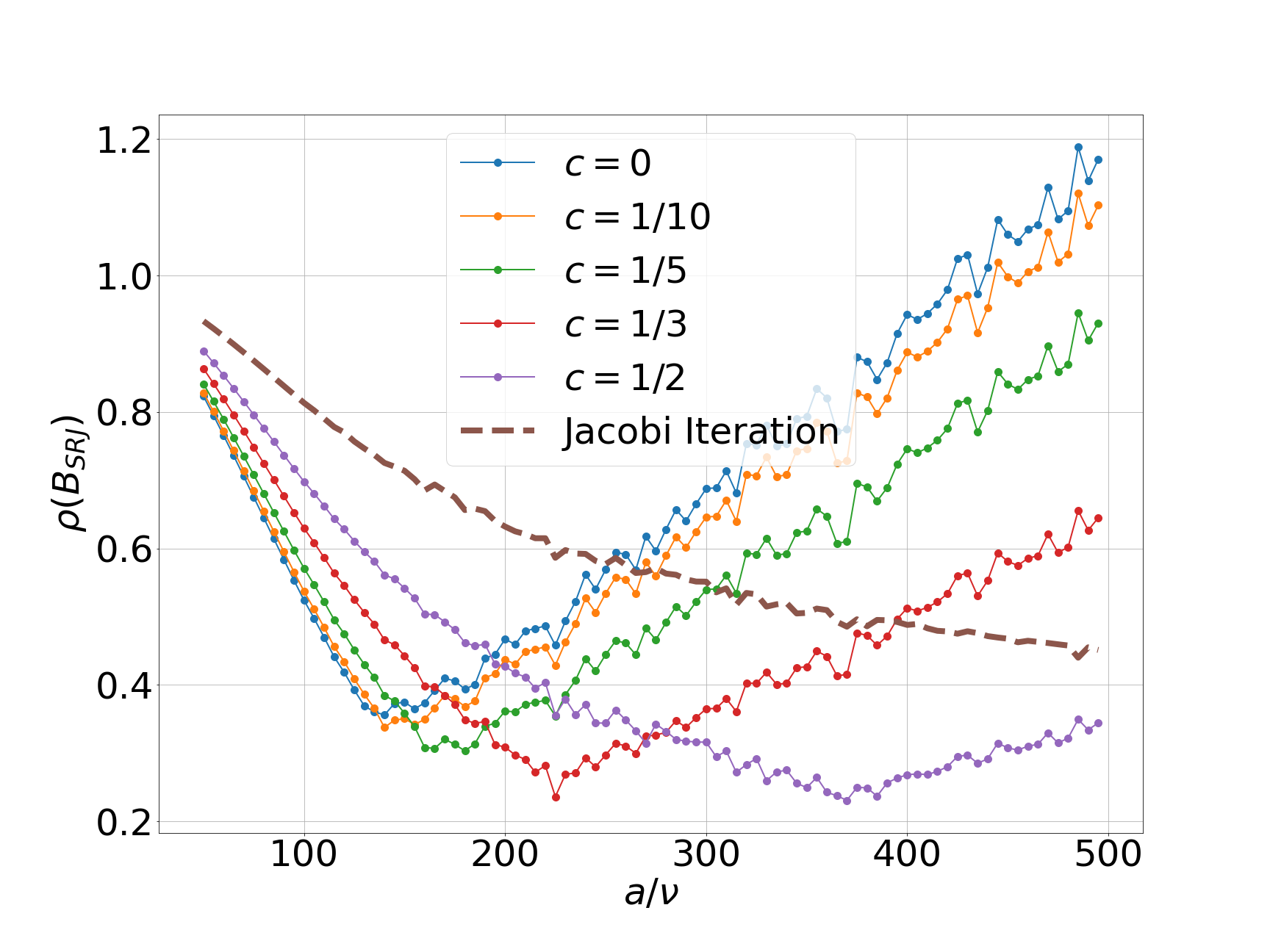}
    \caption{Variation of the spectral radius of the SRJ iteration matrices with advection to diffusion ratio. The schemes here correspond to $M = 5$ relaxation parameters. The spectral radius corresponding to Jacobi iteration (setting the relaxation parameters to 1) is also shown. The scheme with the lowest spectral radius for a given advection to diffusion ratio is expected to provide the fastest convergence.}
    \label{fig:spectral_radius_vs_advection}
\end{figure}

From Figure \ref{fig:spectral_radius_vs_advection}, we observe that the original SRJ scheme corresponding to $c = 0$ (optimized only over the real line) is the best scheme to use for smaller amounts of advection ($a = 50, 100$). As the advection to diffusion ratio increases; however, this scheme provides slower convergence relative to the other schemes. For example, at $\frac{a}{\nu} = 200$, the original $c = 0$ SRJ scheme provides worse convergence relative to all other schemes optimized over a larger elliptical region. This is where we observe the utility of these new SRJ schemes which are optimized for nonsymmetric linear systems. Here the scheme corresponding to $c = 1/3$ provides the fastest convergence. At even larger advection values (e.g. $a = 300$), the original scheme appears to offer no advantage compared to the standard Jacobi iteration. At this point, the scheme corresponding to $c = 1/2$ provides the best convergence. In general, as the amount of advection increases, the schemes which minimize the amplification over a larger elliptical region (larger $c$) exhibit better convergence.

To show how the amount of advection impacts the convergence rate of the different schemes, we plot the spectral radius of the SRJ iteration matrix associated with each scheme, for different advection values $a$. The spectral radius of the iteration matrix $B_{\text{SRJ}}$ determines the asymptotic convergence rate of the SRJ scheme. A smaller spectral radius value is desired for accelerated convergence. We can compute the spectral radius for each SRJ scheme at a given advection value $a$ using Equation \eqref{eqn:srj-convergence}. In particular, given the advection value $a$, the linear system matrix $A$ is known and the Jacobi iteration matrix $B_{\text{J}}$ can be determined. Given the eigenvalues of $B_{\text{J}}$, Equation \eqref{eqn:srj-convergence} can be used to compute the spectral radius of the SRJ iteration matrix for a particular SRJ scheme. This gives a good metric for the anticipated asymptotic convergence behavior, and is also informative of how the schemes will perform relative to each other. 

Figure \ref{fig:spectral_radius_vs_advection} shows the variation of the spectral radius with increasing amounts of advection, for each of the SRJ schemes corresponding to $M = 5$ and $c = 0, 1/10, 1/5, 1/3, 1/2$. The spectral radius associated with the standard Jacobi iterative method (setting all five relaxation factors to one) is also shown. At any given $\frac{a}{\nu}$, the best scheme to use is the one with the smallest spectral radius. Therefore, Figure \ref{fig:spectral_radius_vs_advection} gives a good indication of which scheme is most useful to use for differing advection values. 

For small amounts of advection, the original $c= 0$ SRJ schemes provide the best convergence. When $\frac{a}{\nu}$ exceeds 130, the $c = 0$ scheme no longer has the smallest spectral radius. As a result, other SRJ schemes are expected to provide faster convergence. This is demonstrated by Figure \ref{fig:convergence_a_150}, where we observe that the $c = \frac{1}{10}$ SRJ scheme exhibits the fastest convergence in this case. In general, each scheme has a specific region of advection values for which it provides the fastest convergence. At around $\frac{a}{\nu} = 250$ the spectral radius of the Jacobi iteration is smaller than that of the $c = 0$ SRJ scheme, indicating that standard Jacobi should be faster than this SRJ scheme here. This behavior is verified from Figure \ref{fig:convergence_a_300}, where the $c = 0$ SRJ scheme does no better than Jacobi for $a = 300$. For the highest levels of advection in this plot ($300 < \frac{a}{\nu} < 500$) the SRJ scheme associated with $c = \frac{1}{2}$ is expected to give the best performance. In this high advection regime, we also observe that some of the schemes have an iteration matrix with spectral radius greater than 1, indicating that these scheme will not converge in solving a linear system corresponding to these high advection cases. In particular, the $c = 0, \frac{1}{10}$ are not suitable for convergence for advection-to-diffusion ratio larger than 400. 

The $M = 5$ SRJ schemes derived in this work are shown to provide accelerated convergence for solving the nonsymmetric linear system arising from discretization of the one-dimensional steady advection diffusion PDE, relative to the Jacobi iteration. As the amount of advection grows, the eigenvalue spectrum of the associated Jacobi iteration matrix encompasses a larger region of the complex plane. In these cases, the schemes which minimize the amplification over larger elliptical regions (larger $c$) offer faster convergence compared to the original SRJ schemes which were only derived for the symmetric case. Analysis of the spectral radius of the SRJ iteration matrix corresponding to different schemes also provides valuable insight on which scheme is best to use for a given problem. However, computing the spectral radius of the SRJ iteration matrix may be unfeasible for large scale problems where obtaining the Jacobi iteration matrix eigenvalues is prohibitively expensive. In the next section, we examine a larger two-dimensional problem and show the convergence of the same SRJ schemes in this setting.

\subsection{2D Steady Advection-Diffusion Equation}

\begin{figure}[htbp!]
    \centering
    \includegraphics[width=0.6\textwidth]{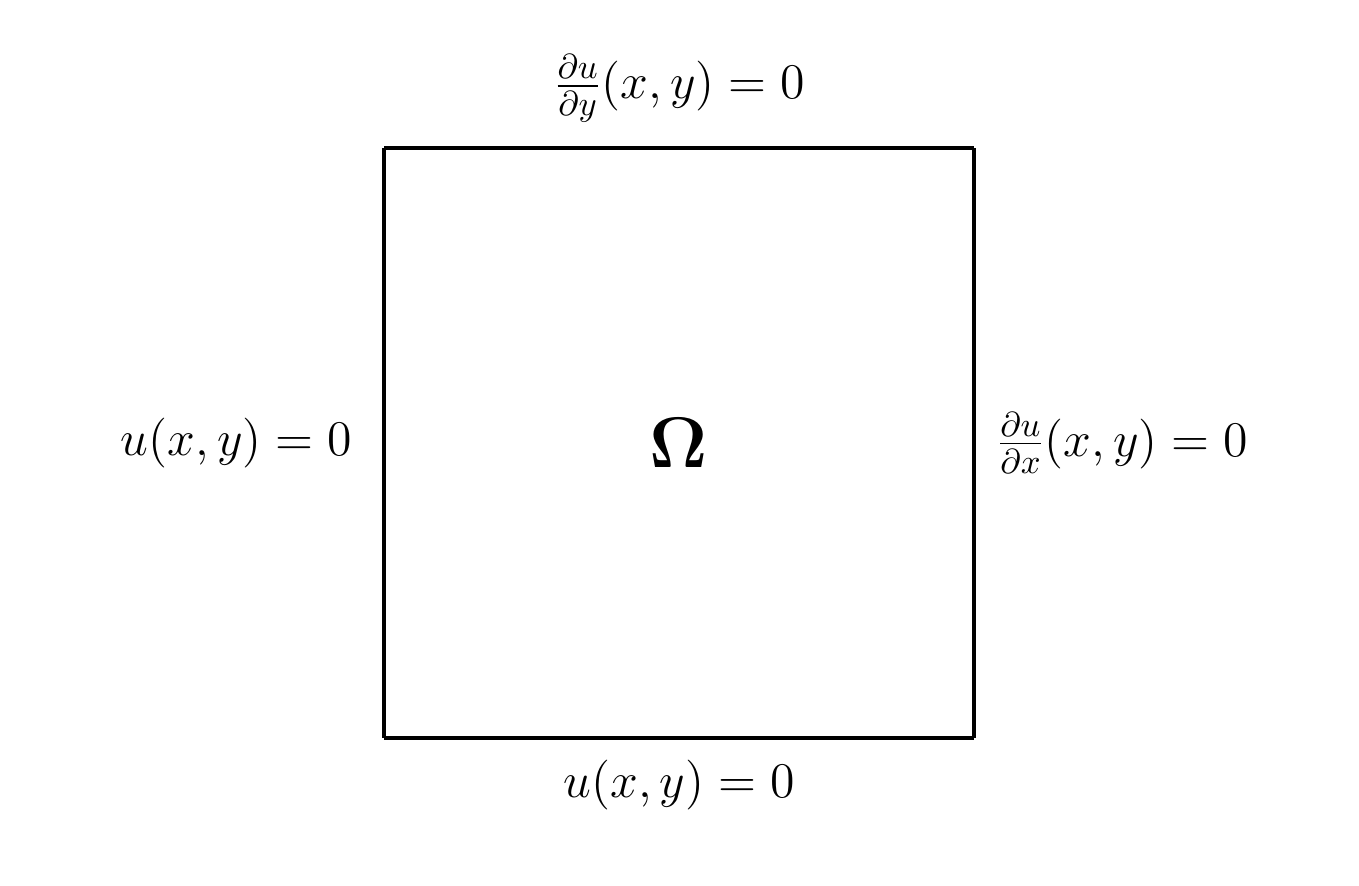}
    \caption{2D Square Domain with homogeneous Dirichlet (bottom and left) and Neumann (top and right) boundary conditions}
    \label{fig:2D-square-domain}
\end{figure}

\begin{figure}[htbp!]
    \begin{subfigure}{.5\textwidth}
        \centering
        \includegraphics[width=\textwidth]{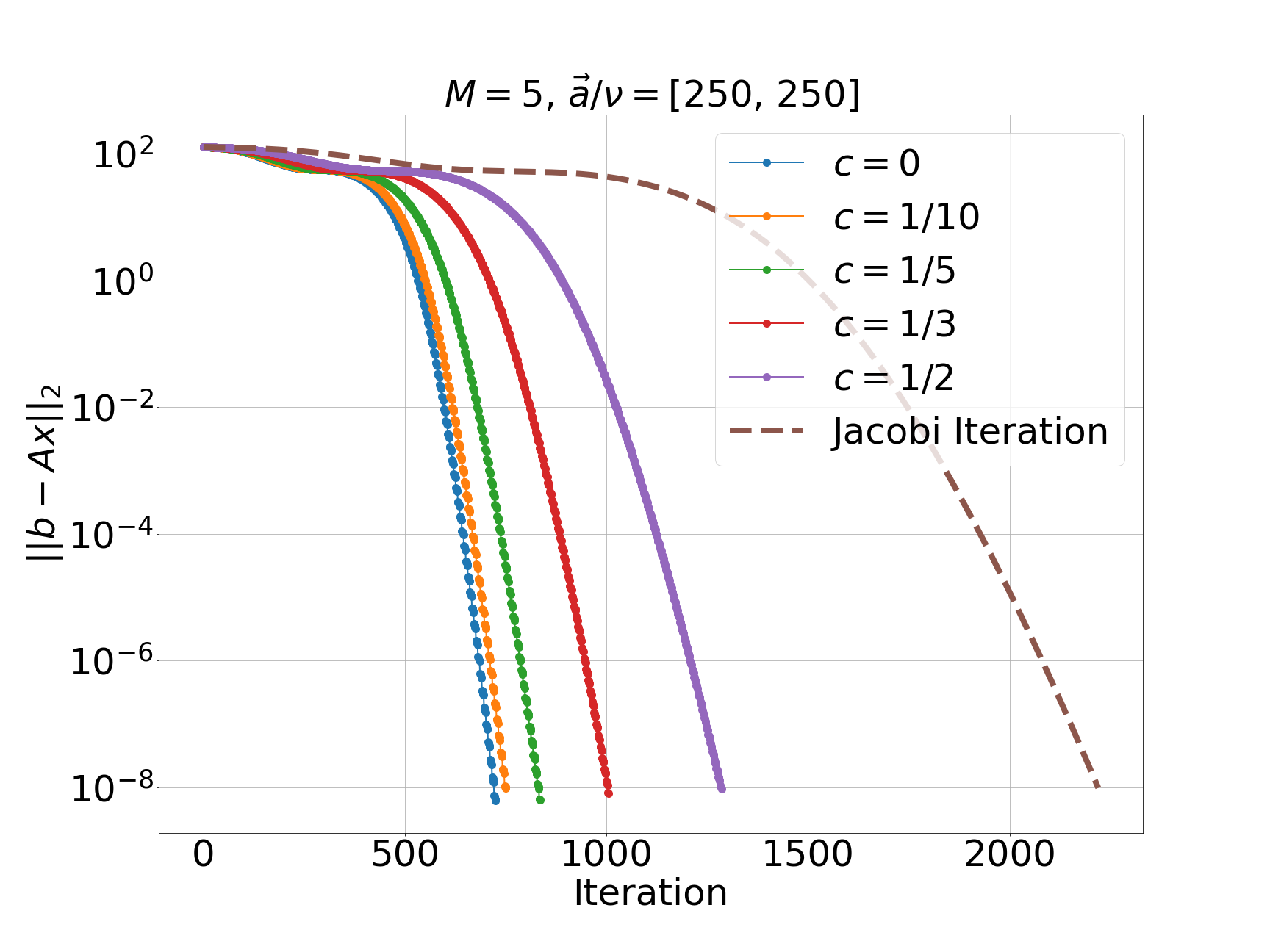}
        \caption{Convergence for $a_{x} = a_{y} = 250$}
        \label{fig:convergence_a_0_2d}
    \end{subfigure}
    \begin{subfigure}{.5\textwidth}
        \centering
        \includegraphics[width=\textwidth]{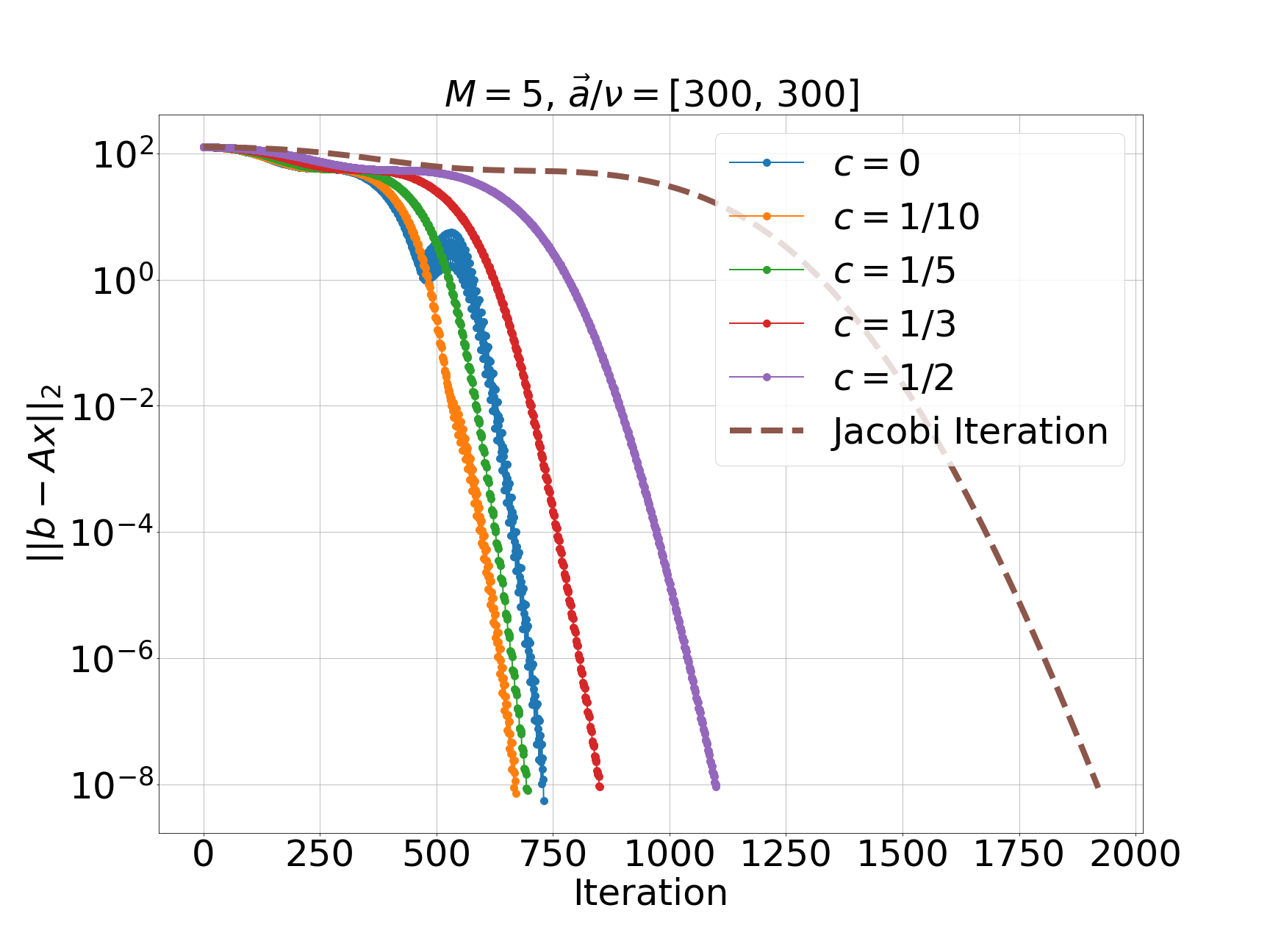}
        \caption{Convergence for $a_{x} = a_{y} = 300$}
        \label{fig:convergence_a_20_2d}
    \end{subfigure}
    \\
    \begin{subfigure}{.5\textwidth}
        \centering
        \includegraphics[width=\textwidth]{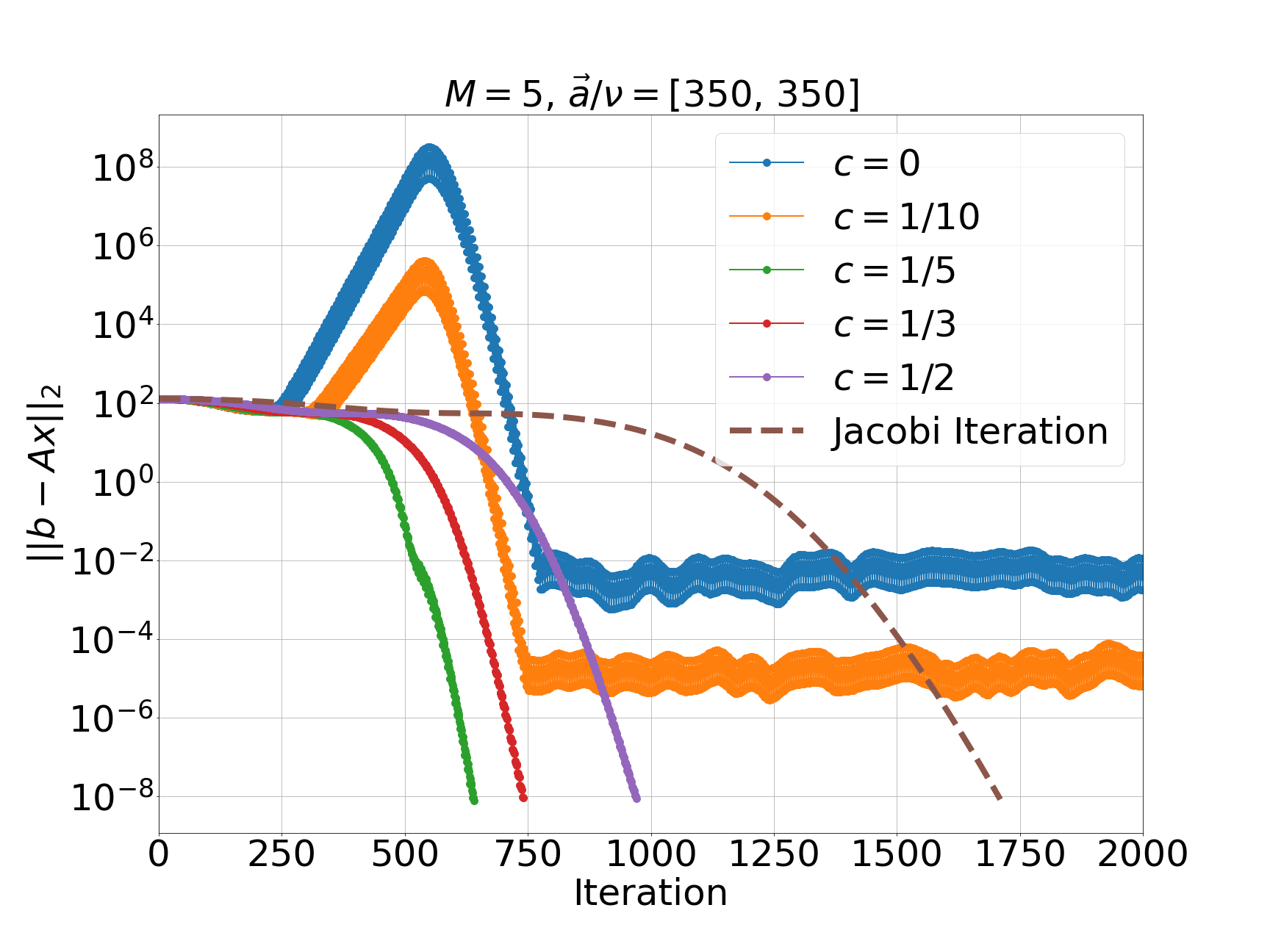}
        \caption{Convergence for $a_{x} = a_{y} = 350$}
        \label{fig:convergence_a_40_2d}
    \end{subfigure}
    \begin{subfigure}{.5\textwidth}
        \centering
        \includegraphics[width=\textwidth]{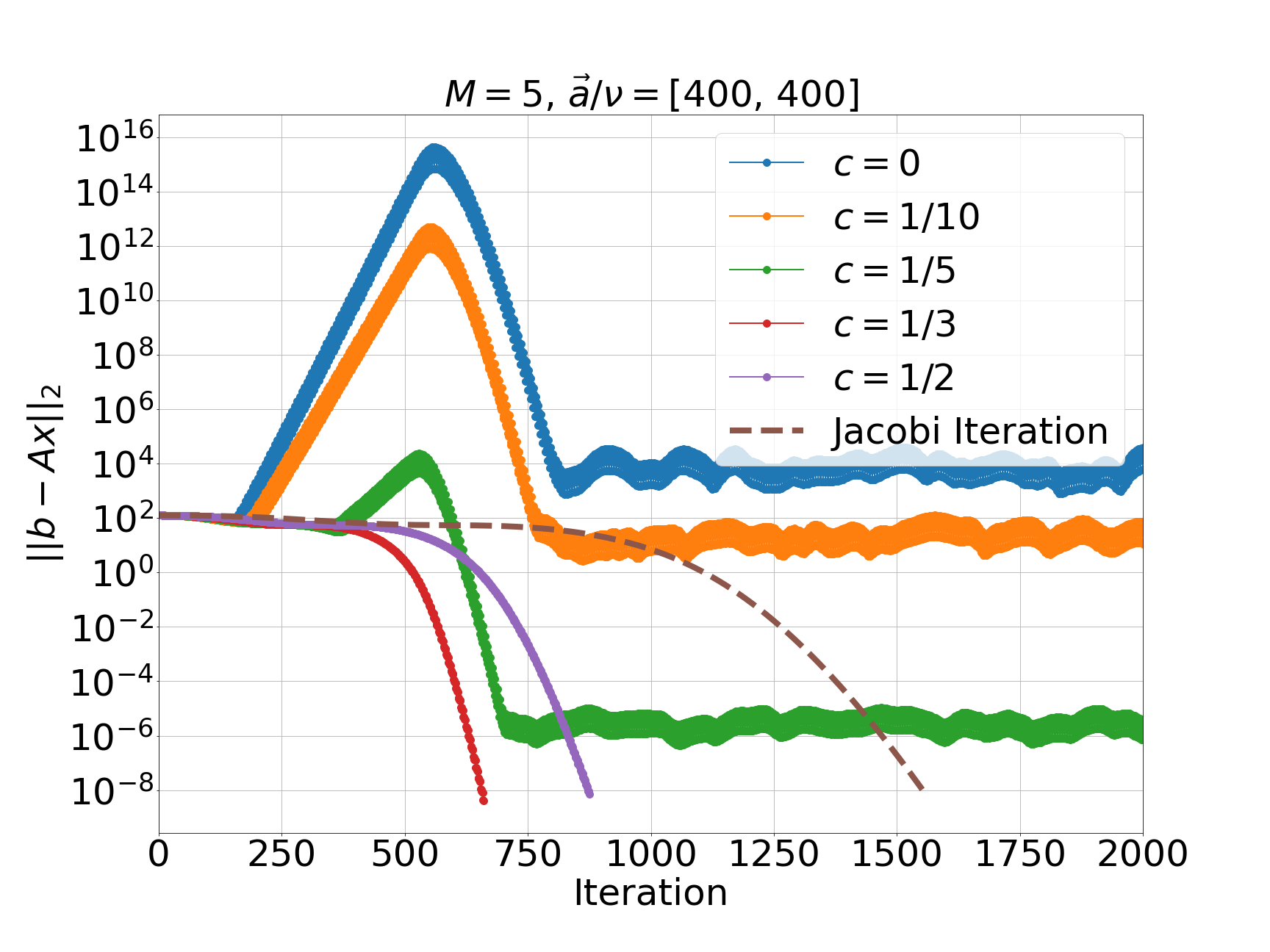}
        \caption{Convergence for $a_{x} = a_{y} = 400$}
        \label{fig:convergence_a_60_2d}
    \end{subfigure}
    \caption{Convergence of the SRJ schemes when solving the 2D advection-diffusion equation with varying amounts of advection. As the amount of advection increases, the schemes optimized over a larger elliptical region exhibit faster convergence relative to the standard Jacobi iteration and SRJ optimized only over the real axis ($c = 0$). Several of the schemes stagnate for cases of high advection. In all cases, one of the SRJ schemes outperforms the Jacobi iteration method.}
    \label{fig:convergence_of_advection_diffusion_2D}
\end{figure}

We consider the two dimensional steady advection-diffusion equation given by Equation \eqref{eqn:2d-advection-diffusion}
\begin{equation}
    - \nu \nabla^2 u(x) + \vec{a} \cdot \nabla u(x) = f(x)
    \label{eqn:2d-advection-diffusion}
\end{equation}
where $\nu$ is the diffusion coefficient and $\vec{a} = [a_{x}, a_{y}]^T$ is a vector containing the advection coefficients in the $x$ and $y$ directions. We employ a finite difference discretization in a unit square domain to discretize the PDE in Equation \eqref{eqn:2d-advection-diffusion} into a set of linear equations. We assume homogenous Dirichlet boundary conditions at the bottom and left side, and homogenous Neumann boundary conditions at the right and top side as given by Equations \eqref{eqn:2D-bc-dirichlet}-\eqref{eqn:2D-bc-neumann}.
\begin{align}
    & \text{Dirichlet BCs}: u(x = 0, y) = 0, \ u(x, y = 0) = 0, 
    \label{eqn:2D-bc-dirichlet} \\
    & \text{Neumann BCs}: \frac{\partial u(x,y)}{\partial x}\bigg|_{x = 1, y} = 0, \ \frac{\partial u(x,y)}{\partial y}\bigg|_{x, y = 1} = 0, \
    \label{eqn:2D-bc-neumann}
\end{align}
A schematic of the domain and boundary conditions is shown in Figure \ref{fig:2D-square-domain}.
We consider a two-dimensional finite difference stencil for discretizing the continuous derivatives in the PDE in Equation \eqref{eqn:2d-advection-diffusion}. In particular, a second order central difference discretization is employed for the second derivatives 
\begin{align}
    &\frac{d^{2} u}{dx^2} \approx \frac{u_{i+1,j} - 2u_{i,j} + u_{i-1,j}}{\Delta x^2} \\
    &\frac{d^{2} u}{dy^2} \approx \frac{u_{i,j+1} - 2u_{i,j} + u_{i,j-1}}{\Delta y^2} 
\end{align}
while an upwinding scheme is used for the first derivative in both directions as follows:
\begin{align}
    &\frac{du}{dx} \approx
    \begin{cases}
        \frac{u_{i,j} - u_{i-1.j}}{\Delta x} \ \text{if} \ a_{x} > 0 \\
        \frac{u_{i+1,j} - u_{i,j}}{\Delta x} \ \text{if} \ a_{x} < 0 \\
    \end{cases} \\
    &\frac{du}{dy} \approx
    \begin{cases}
        \frac{u_{i,j} - u_{i,j-1}}{\Delta y} \ \text{if} \ a_{y} > 0 \\
        \frac{u_{i,j+1} - u_{i,j}}{\Delta x} \ \text{if} \ a_{y} < 0
    \end{cases}
\end{align}
where $u_{i,j}$ represents the solution at the $(i,j)$ grid point in our two dimensional grid, and $\Delta x$ and $\Delta y$ are the grid spacing in the $x$ and $y$ directions. The finite difference discretization leads to a pentadiagonal system of equations $Ax = b$ where the matrix $A$ is given by 
\[
A = 
\begin{pmatrix}
a_{0} & a_{1} & & a_{N} & & \\
a_{-1} & a_{0} & a_{1} & & \ddots & \\
 & \ddots & \ddots & \ddots & & a_{N} \\
a_{-N} &  & \ddots & \ddots & \ddots \\
 & \ddots &  & \ddots & \ddots & a_{1} \\
  &  &  a_{-N} & & a_{-N} & a_{0} 
\end{pmatrix} \ \
\]
where (assuming positive $a_{x}$ and $a_{y}$)
\begin{align*}
    a_{-N} &= -\frac{\nu}{\Delta y^2} - \frac{a_{y}}{\Delta y} \\
    a_{-1} &= -\frac{\nu}{\Delta x^2} - \frac{a_{x}}{\Delta x} \\
    a_{0} &= \frac{2\nu}{\Delta x^2} + \frac{2\nu}{\Delta y^2} + \frac{a_{x}}{\Delta x} + \frac{a_{y}}{\Delta y} \\
    a_{1} &= -\frac{\nu}{\Delta x^2} \\
    a_{N} &= -\frac{\nu}{\Delta y^2} 
\end{align*}
except for rows corresponding to boundary conditions which must be adjusted. We discretize the two dimensional steady advection-diffusion equation with $N_{x} = 256$ by $N_{y} = 256$ interior DOFs. We explore different values of the advection coefficient and solve the resulting linear system using the $M = 5$ SRJ schemes. In particular, we set $a_{x} = a_{y} = 250, 300, 350, 400$ and $\nu = 1$ and examine the convergence behavior of each of the schemes for solving these linear systems. The convergence results are shown in Figure \ref{fig:convergence_of_advection_diffusion_2D}, showing the $L_{2}$ residual norm at each SRJ iteration until a residual below a tolerance value of 1E-8 is achieved.

For the smallest advection value of $a_{x} = a_{y} = 250$ explored in this study, the original $c = 0$ SRJ scheme provides the fastest convergence. The schemes associated with larger $c$ (minimizing amplification over larger ellipse regions) demonstrate convergence behavior which is progressively worse as $c$ increases. However, all schemes outperform the standard Jacobi iteration here. When $a_{x} = a_{y} = 300$, the $c = \frac{1}{10}, \frac{1}{5}$ schemes provide faster convergence than the $c = 0$ scheme, illustrating that the linear system now has complex eigenvalues which are more well-suited for these new schemes. At $a_{x} = a_{y} = 350$, the $c = 0, \frac{1}{10}$ schemes begin to diverge and eventually stagnate without achieving our desired convergence. The $c = \frac{1}{5}$ scheme provides the best convergence here. This indicates that the Jacobi iteration matrix eigenvalues are distributed such that they lie outside the elliptical region associated with the $c = \frac{1}{10}$ ellipse but inside the region corresponding to $c = \frac{1}{5}$. At $a_{x} = a_{y} = 400$, the schemes associated with $c = 0, \frac{1}{10}, \frac{1}{5}$ all diverge, then start to converge, and eventually stagnate and are unable to achieve the desired tolerance. The $c = \frac{1}{3}$ scheme demonstrates the best convergence here, indicating that the iteration matrix eigenvalue distribution has expanded in the complex plane. In all cases explored, there is an SRJ scheme which provides accelerated convergence relative to the Jacobi iteration scheme.

Our two dimensional advection diffusion results illustrate the utility of the newly developed SRJ schemes. In cases of low advection, the original $c = 0$ SRJ schemes usually provide the fastest convergence behavior. However, in cases of high advection, this scheme tends to show divergence behavior and eventually stagnates, being unable to achieve the desired tolerance. In these cases, the other schemes are able to achieve good performance and still outperform the standard Jacobi iteration and its convergence behavior. With more advection, the schemes associated with larger $c$ do a better job of achieving faster convergence.


\section{Conclusion}

The main goal of this paper is to broaden the applicability of the Scheduled Relaxation Jacobi (SRJ) method which, until this point, has mainly been restricted to solving partial differential equations which are elliptic. In this work, we have developed a methodology for constructing SRJ schemes for solving non-elliptic PDEs, or equivalently, the nonsymmetric linear systems arising from their discretization via numerical solution to an optimization problem. The schemes developed in this work are provided in Appendix \ref{sec:all-srj-schemes}, and are an extension of SRJ schemes previously developed for solving elliptic PDEs \cite{islam2020}. The schemes are obtained by solving a constrained minimization problem which seeks to minimize the maximum amplification within elliptical regions of the complex plane. The schemes developed in this work are used to solve the steady advection-diffusion equation discretized by finite difference schemes in 1D and 2D, which grow progressively more nonsymmetric as the ratio of advection to diffusion is increased. As the amount of advection increases, using the modified schemes greatly improves the convergence of the linear system compared to the standard Jacobi iteration and the original SRJ schemes derived for symmetric problems. While we have provided several schemes for specified number of relaxation parameters $M$ and ellipse regions corresponding to $c$, other schemes could be derived by adjusting these values and solving the optimization problem in these cases. Practitioners can directly utilize the SRJ scheme parameters given in Appendix \ref{sec:all-srj-schemes} for their own work, or derive new schemes for different $M$ and $c$ which are more suitable for their applications. The optimization framework shown here may even be used to develop schemes suitable for problems whose Jacobi iteration matrix eigenvalue spectra resemble non-ellipse regions. 

The SRJ method shows great promise as a high performance linear solver method due to its potential for massive parallelization on high performance computing architectures. Extending the available SRJ schemes to the nonsymmetric case makes the method suitable for solving more general problems arising from a wider array of applications which are of interest to scientists and engineers. 

\clearpage
\printbibliography

\clearpage
\appendix

\section{Derivation of Constraint Jacobian for Optimization Routine}
\label{sec:appendix1-jacobian-derivation}

This section details the derivation for the constraint Jacobian. We begin with the amplification polynomial $G_{M}$ which is given by
\begin{equation}
    G_{M} (\lambda) = \prod_{i=1}^{N} (1-\omega_i) + \omega_i \lambda
\end{equation}
The inputs $\lambda$ to the amplification polynomial are permitted to be complex numbers (i.e. $\lambda = x + iy$). Therefore, the polynomial evaluated at some complex $\lambda$ is the product of many complex values with real and imaginary part $A_{i}$ and $B_{i}$ as follows
\begin{equation} 
    G_{M} (\lambda) = \prod_{i=1}^{N} \underbrace{(1-\omega_i) + \omega_i x}_{A_i} + i \underbrace{\omega_i y}_{B_i}
    \label{eqn:amplification-complex}
\end{equation}
The constraint for our optimization is enforced at a collection of test points (see Figure \ref{fig:test_points}) and is given by
\begin{equation}
    c (\vec{x}) \equiv \bar{g}^2 - |G_{M}(\lambda;w_j)|^2 
\end{equation}
where $c (\vec{x}) \geq 0$, and the solution vector $\vec{x} = [w_j; \bar{g}]^{T}$ consists of the $M$ relaxation factors and bounding value $\bar{g}$. The Jacobian vector associated with this constraint function consists of $M+1$ entries, each corresponding to the derivative of the constraint function with respect to the optimization parameters in $\vec{x}$
\begin{equation}
    J = \frac{\partial c(\vec{x})}{\partial \vec{x}} = \bigg[\underbrace{\frac{\partial c(\vec{x})}{\partial \omega_j}}_{M \ \text{terms}}, \underbrace{\frac{\partial c(\vec{x})}{\partial \bar{g}}}_{1 \ \text{term}}\bigg]^{T}
    \label{eqn:jacobian-vector}
\end{equation}
The remainder of the analysis is primarily focused on deriving the $M$ terms corresponding to $\frac{\partial c(\vec{x})}{\partial \omega_j}$. For brevity, we will not explicity write down function inputs (i.e $c({\vec{x}}) \equiv c$, $G_{M}(\lambda;w_j) \equiv G_{M}$). Given that $G_{M}$ can take on complex values, we can write its squared magnitude as the sum of squares of its real and imaginary components as follows
\begin{equation}
    |G_{M}|^2  = \text{Re}^2[G_{M}] + \text{Im}^2[G_{M}]
\end{equation}
Therefore, the constraint function can be written as
\begin{equation}
    c \equiv \bar{g}^2 - \text{Re}^2[G_{M}] - \text{Im}^2[G_{M}]
    \label{eqn:constraint-2}
\end{equation}
The derivative of the constraint with respect to the relaxation parameters $\omega_j$ is
\begin{equation}
    \frac{\partial c}{\partial \omega_j} = -2 \ \text{Re}[G_{M}] \frac{\partial (\text{Re}[G_{M}])}{\partial \omega_j} - 2 \ \text{Im}[G_{M}] \frac{\partial (\text{Im}[G_{M}])}{\partial \omega_j}
    \label{eqn:jacobian-terms}
\end{equation}
We determine expressions for the real and imaginary parts of the function $G_{M}$, and can afterwards obtain the derivative of these expressions with respect to relaxation factors $w_j$, providing the ingredients to compute the Jacobian. Starting from Equation \eqref{eqn:amplification-complex}, we can express the function $G_{M}$ has a product of complex exponentials as follows
\begin{align}
    G_{M} = \prod_{i=1}^{N} A_i + i B_i = \prod_{i=1}^{N} C_i e^{i \theta_i}
\end{align}
where 
\begin{align}
    \label{eqn:ci}
    C_i &= \sqrt{A_i^2 + B_i^2} \\
    \label{eqn:thetai}
    \theta_i &= \tan^{-1} \left(\frac{B_i}{A_i}\right) 
\end{align}
We have used the fact that any complex number can be written as a complex exponential. Furthermore, we can combine the product of complex exponential into a single complex exponential
\begin{align}
    G_{M} = \prod_{i=1}^{N} C_i e^{i \theta_i} = D^{*} e^{i \theta^{*}}
\end{align}
where 
\begin{align}
    \label{eqn:dstar}
    D^{*} =  \prod_{i=1}^{N} C_i \\
    \label{eqn:thetastar}
    \theta^{*} = \sum_{i=1}^{N} \theta_i
\end{align}
The real and imaginary parts of $G_{M}$ can now be easily obtained. Euler's formula states that $e^{i \theta} = \cos{\theta} + i \sin{\theta}$. Therefore, it follows that $\text{Re}[e^{i \theta}] = \cos{\theta}$ and $\text{Im}[e^{i \theta}] = \sin{\theta}$. Following this logic, we can express the real and imaginary parts of $G_{M}$ as
\begin{align}
    \label{eqn:re_g}
    \text{Re}[G_{M}] = D^{*} \cos(\theta^{*}) \\
    \label{eqn:im_g}
    \text{Im}[G_{M}] = D^{*} \sin(\theta^{*})
\end{align}
We can now obtain the derivatives of these expressions with respect to the relaxation factors $w_j$. Since $D^{*}$ and $\theta^{*}$ both depend on the relaxation factors $\omega_i$ (due to their dependence on $C_i$, $\theta_i$ which depend on $A_i$, $B_i$), we can obtain the derivative using the product rule and chain rule which results in
\begin{align}
    \label{eqn:re_g_deriv}
    \frac{\partial (\text{Re}[G_{M}])}{\partial \omega_j}  &= \frac{\partial D^{*}}{\partial \omega_j} \cos(\theta^{*}) - D^{*} \sin(\theta^{*}) \frac{\partial \theta^{*}}{\partial \omega_j}  \\
    \frac{\partial (\text{Im}[G_{M}])}{\partial \omega_j} &= \frac{\partial D^{*}}{\partial \omega_j} \sin(\theta^{*}) + D^{*} \cos(\theta^{*}) \frac{\partial \theta^{*}}{\partial \omega_j} 
    \label{eqn:im_g_deriv}
\end{align}
To compute the derivatives above, we must compute the derivative of $D^{*}$ and $\theta^{*}$ with respect to $\omega_j$. Starting with $D^{*}$ and using Equation \eqref{eqn:dstar}
\begin{equation}
    \frac{\partial D^{*}}{\partial \omega_j} = \frac{\partial}{\partial \omega_j} \left[ \prod_{i=1}^{N} C_i \right] = \frac{\partial C_j}{\partial \omega_j} \left[ \prod_{i=1, i \neq j}^{N} C_i \right]
\end{equation}
Note that only one of the $C_i$ depend on a particular $\omega_j$, when the two indices are equal. We can determine $\frac{\partial C_j}{\partial \omega_j}$ as follows
\begin{align}
    \frac{\partial C_j}{\partial \omega_j} &= \frac{\partial}{\partial \omega_j} \left[\sqrt{A_j^2 + B_j^2} \right] \\
    &= \frac{1}{2}  \left( A_j^2 + B_j^2 \right)^{-1/2} \left[ 2A_j \frac{dA_j}{d\omega_j} + 2B_j \frac{dB_j}{d\omega_j} \right] \\
    &= \frac{A_j (x-1) + B_j y}{\sqrt{A_j^2 + B_j^2}}
\end{align}
where we have used the fact that $A_j = 1-\omega_j+\omega_j x$ and $B_j = \omega_j y$ so that $\frac{dA_j}{d\omega_j} = x-1$ and $\frac{dB_j}{d\omega_j} = y$, where $x = \text{Re}(\lambda)$, $y = \text{Im}(\lambda)$. The final expression for the derivative of $D^{*}$ with respect to $\omega_j$ is 
\begin{equation}
    \frac{\partial D^{*}}{\partial \omega_j} = \frac{A_j (x-1) + B_j y}{\sqrt{A_j^2 + B_j^2}}  \left[ \prod_{i=1, i \neq j}^{N} C_i \right]
    \label{eqn:dDstar_dwj}
\end{equation}
We now determine the derivative of $\theta^{*}$ with respect to $\omega_j$. 
\begin{align}
    \frac{\partial \theta^{*}}{\partial \omega_j} = \frac{\partial }{\partial \omega_j} \left[\sum_{i=1}^{N} \theta_i \right] = \frac{\partial \theta_j}{\partial \omega_j} &= \frac{\partial }{\partial \omega_j} \left[\tan^{-1} \left(\frac{B_j}{A_j} \right) \right] \\
    &= \frac{1}{1 + \left(\frac{B_j^2}{A_j^2} \right)} \frac{\frac{\partial B_j}{\partial w_j} A_j + B_j \frac{\partial A_j}{\partial w_j}}{A_j^2} \\
    & = \frac{A_j y + B_j (x-1)}{A_j^2 + B_j^2}
    \label{eqn:dthetastar_dwj}
\end{align}
In the above derivation, it is true that only one $\theta_i$ term depends on a particular $\omega_j$, when the two indices are equal. Furthermore, we employed the following derivative 
\begin{equation}
    \frac{d}{dx} \left[ \tan^{-1} x \right] = \frac{1}{1+x^2}
\end{equation}
and used the chain rule to complete the derivation. Given Equations \eqref{eqn:dDstar_dwj} and \eqref{eqn:dthetastar_dwj}, we can now evaluate the derivatives of $\text{Re}[G_{m}]$ and $\text{Im}[G_{m}]$ with respect to $\omega_j$ given in Equations \eqref{eqn:re_g_deriv}-\eqref{eqn:im_g_deriv}.

We now have all of the ingredients necessary to compute the $M$ entries of the Jacobian vector corresponding to $\frac{\partial c}{\partial \omega_j}$. The full procedure to compute these entries is as follows. Given the current scheme parameters $\omega_i$ and some potentially complex input $\lambda$ at which to compute the Jacobian entries (for our optimization problem these inputs correspond to a test point location at which to enforce the constraint) the entries of the Jacobian vector $\frac{\partial c}{\partial \omega_j}$ can be obtained following the steps below
\begin{enumerate}
    \item Determine $x = \text{Re}(\lambda)$ and $y = \text{Im}(\lambda)$ and compute the coefficients $A_i = (1 - \omega_i) + \omega_i x$, $B_i = \omega_i y$
    \item Compute $C_i, \theta_i$ using Equations \eqref{eqn:ci}, \eqref{eqn:thetai} and $D^{*}, \theta^{*}$ using Equations \eqref{eqn:dstar}, \eqref{eqn:thetastar}
    \item Compute $\frac{\partial D^{*}}{\partial \omega_j}$ and $\frac{\partial \theta^{*}}{\partial \omega_j}$ using Equations \eqref{eqn:dDstar_dwj} and \eqref{eqn:dthetastar_dwj}
    \item Compute $\text{Re}[G_{M}]$, $\text{Im}[G_{M}]$ using Equations \eqref{eqn:re_g}, \eqref{eqn:im_g}, and $\frac{\partial \text{Re}[G_{M}]}{\partial \omega_j}$, $\frac{\partial \text{Im}[G_{M}]}{\partial \omega_j}$ using Equations \eqref{eqn:re_g_deriv}, \eqref{eqn:im_g_deriv}
    \item Compute $\frac{\partial c}{\partial \omega_j}$ using Equation \eqref{eqn:jacobian-terms}
\end{enumerate}
We can now compute all terms associated with $\frac{\partial c}{\partial \omega_j}$ in the Jacobian vector in Equation \eqref{eqn:jacobian-vector}. The final entry of the Jacobian vector is given by 
\begin{equation}
    \frac{\partial c}{\partial \bar{g}} = 2\bar{g}
\end{equation}
This fully specifies the constraint Jacobian. The constraint and constraint Jacobian can be determined at each test point location for the optimization routine.

\newpage
\section{Comparison of SRJ schemes for solving stiff systems}
\label{sec:scheme-slope}

The ability of a SRJ scheme to converge for a stiff problem can be characterized by the slope of its corresponding amplification polynomial at $\lambda = 1$. If the slope is higher, the spectral radius of the SRJ scheme for the stiff problem is likely to deviate further away from 1, corresponding to a faster asymptotic convergence rate. Table \ref{tab:slope_table} shows this slope for all schemes associated with $M = 2$ to $M = 20$ relaxation factors and $c = 0, \frac{1}{10}, \frac{1}{5}, \frac{1}{3}, \frac{1}{2}$. 

Schemes associated with larger $M$ are better able to handle stiff problems. An explanation for this is that the amplification polynomials associated with larger $M$ bound a larger portion of the $\lambda \in (-1,1)$ region. A stiff problem would have an iteration matrix eigenvalue close to 1, so schemes associated with larger $M$ are better able to "capture" this eigenvalue and bound its associated amplification. Schemes associated with larger $c$ are not able to handle stiff systems as well. This is because these schemes prioritize bounding the amplification in a wider region of the complex plane, rather than a wider region of the real axis. The spread of the iteration matrix eigenvalues is important for determining which SRJ scheme should be used.

For another point of comparison, the amplification corresponding to the Jacobi scheme with $M$ relaxation factors (all set to one) has a slope of $M$ at $\lambda = 1$. This can be shown by taking the $M$ order amplification polynomial (Equation \eqref{eqn:amplification-M-scheme}) and setting all relaxation factors to one, which results in an amplification of $G_{M}(\lambda) = \lambda^M$ for $M$ steps of Jacobi iteration. For any SRJ scheme in Table \ref{tab:slope_table} corresponding to a particular $M$ and $c$, the slope of the scheme is always greater than $M$. This indicates that all schemes developed in this work accelerate the convergence of standard Jacobi iteration when solving stiff linear systems. 

\vspace{0.2in}

\begin{table}[htbp!]
\caption{Measure of SRJ scheme's ability to handle stiff linear system. Schemes associated with higher $M$ and smaller $c$ are best for solving stiff systems.}
\centering
\begin{tabular}{ |c|c|c|c|c|c|c| } 
\hline
\diagbox{$M$}{$c$} & 0 & $\frac{1}{10}$ & $\frac{1}{5}$ & $\frac{1}{3}$ & $\frac{1}{2}$ & \text{Jacobi} \\
\hline
2 &  2.276 &  2.269 & 2.246 & 2.195 & 2.101 & 2 \\
\hline
3 & 4.951 & 4.905 & 4.770 & 4.485 & 4.035 & 3 \\
\hline
4 & 8.696 & 8.539 & 8.109 & 7.283 & 6.168 & 4  \\
\hline
5 & 13.510 & 13.121 & 12.112 & 10.371 & 8.349 & 5 \\
\hline
6 & 19.393 & 18.588 & 16.624 & 13.598 & 10.521 & 6  \\
\hline
7 & 26.346 & 24.871 & 21.509 & 16.877 & 12.671 & 7  \\
\hline
8 & 34.369 & 31.894 & 26.652 & 20.161 & 14.798 & 8 \\
\hline
9 & 43.461 & 39.582 & 31.962 & 23.429 & 16.906 & 9  \\
\hline
10 & 53.624 & 47.855 & 37.373 & 26.674 & 18.998 & 10 \\
\hline
11 & 64.855 & 56.640 & 42.837 & 29.896 & 21.077 & 11 \\
\hline
12 & 77.156 & 65.866 & 48.323 & 33.094 & 23.145 & 12 \\
\hline
13 & 90.528 & 75.465 & 53.809 & 36.271 & 25.204 & 13 \\
\hline
14 & 104.968 & 85.380 & 59.281 & 39.431 & 27.256 & 14 \\
\hline
15 & 120.479 & 95.552 & 64.735 & 42.574 & 29.302 & 15 \\
\hline
16 & 137.059 & 105.933 & 70.164 & 45.704 & 31.342 & 16 \\
\hline
17 & 154.709 & 116.487 & 75.570 & 48.821 & 33.379 & 17 \\
\hline
18 & 173.428 & 127.172 & 80.951 & 51.928 & 35.411 & 18 \\
\hline
19 & 193.217 & 137.958 & 86.307 & 55.025 & 37.441 & 19 \\
\hline
20 & 214.079 & 148.821 & 91.640 & 58.114 & 39.468 & 20 \\
\hline
\end{tabular}
\label{tab:slope_table}
\end{table}

\newpage
\section{SRJ schemes for nonsymmetric matrices}
\label{sec:all-srj-schemes}

This appendix lists all of the SRJ schemes developed in this work. The SRJ schemes associated with $M = 2$ up to $M = 20$ parameters, for each of $c = 0, \frac{1}{10}, \frac{1}{5}, \frac{1}{3}, \frac{1}{2}$, are listed in Tables \ref{tab:srj-schemes-c0}-\ref{tab:srj-schemes-c1_2}.

\begin{table}[htbp!]
\caption{SRJ scheme parameters for $c = 0$}
\begin{center}
\begin{tabular}{ |c|p{150mm}| } 
\hline
$M$ & \textbf{SRJ scheme parameters} \\
\hline
2 & 0.5690356, 1.70710677\\
\hline
3 & 3.49402001, 0.53277775, 0.9245737 \\
\hline
4 & 0.70836006, 1.46555904, 6.00294106, 0.51881773 \\
\hline
5 & 2.17132943, 0.97045898, 0.51215172, 9.23070087, 0.62486987 \\
\hline
6 & 0.78456896, 3.03760966, 13.17650299, 0.50847872, 1.30215429, 0.58365629  \\
\hline
7 & 0.50624677, 0.98455485, 1.69891723, 0.69311373, 17.84007874, 4.063045, 0.56014437 \\
\hline
8 & 5.24709171, 23.2213142, 0.83002589, 0.54540533, 0.50479132, 0.64064728, 1.22055571, 2.15908451 \\
\hline
9 & 0.60746266, 0.50379038, 1.49082726, 6.58949616, 0.99056048, 0.74168003, 0.53553197, 2.68192018, 29.3201555 \\
\hline
10 & 0.50307289, 0.85987944, 36.13657316, 0.58501073, 0.68577217, 1.17301688, 3.26705563, 8.09012714, 1.79453894, 0.52858345 \\
\hline
11 & 9.74878298, 3.91423663, 43.67026622, 0.52350202, 0.77740428, 0.99363986, 0.64784306, 0.56904981, 1.37652217, 2.13122571, 0.50254084 \\
\hline
12 & 0.51967186, 0.502136, 1.60062418, 0.88089529, 1.14212703, 0.5572681, 11.56563689, 0.7214957, 0.62078228, 4.62343187, 2.50067832, 51.92170993 \\
\hline
13 & 5.39452922, 60.89068271, 13.54056493, 0.51671135, 0.68157772, 0.54830598, 0.80462771, 0.99542885, 1.84503098, 0.60071874, 2.90272506, 0.50182066, 1.30484854 \\
\hline
14 & 0.51437483, 0.50157028, 1.12052206, 1.4815068, 0.75022361, 70.57717917, 3.33726813, 0.89647097, 0.58538736,  0.5413206, 6.22748084, 15.67354847, 2.10956602, 0.65193728 \\
\hline
15 & 0.50136817, 80.98119527, 0.82600199, 1.25587196, 0.62924226, 2.39411386, 17.96457474, 3.80424262, 0.7098291, 1.6719108, 0.57338285, 7.12225455, 0.99655627, 0.51249792, 0.53576477 \\
\hline
16 & 92.10272776, 20.4136347, 4.30360397, 0.67888335, 0.53126985, 8.07882802, 0.61143075, 0.56379208, 0.90846298,  1.10458234, 0.50120269, 0.51096714, 1.87593302, 1.40127987, 2.69859625, 0.77372736 \\
\hline
17 & 0.50106549, 0.65457349, 4.83532073, 103.94177396, 1.22034776, 0.55599891, 0.99731168, 0.5275797, 1.5566107, 0.59716511, 0.8432046, 0.5097021, 3.02295855, 23.02072172, 0.73369365, 2.0934855, 9.09718533 \\
\hline
18 & 0.50095048, 0.50864452, 116.49833154, 0.54957426, 1.34371329, 3.36716161, 0.70224149, 25.78583082, 1.72176974,   0.52451158, 0.58554342, 0.91797527, 1.09234748, 2.32450583, 0.63507926, 0.79326944, 10.17731485, 5.3993702 \\
\hline
19 & 3.73117675, 5.99573557, 0.75405501, 129.77239844, 0.99784208, 0.5442113, 11.31920792, 0.50085313, 2.56894868,   1.47457999, 0.50775129, 0.6191755, 0.85733643, 1.19342964, 1.89668891, 0.52193218, 28.70895814, 0.57593767, 0.67700928 \\
\hline 
20 & 0.50077028, 0.5069904, 1.30047524, 12.52311806, 1.6129012, 143.76536998, 1.08267819, 31.7907227, 4.11506368, 0.519743, 0.53968667, 6.62454016, 2.08135409, 2.82683304, 0.92571219, 0.72260429, 0.56790042, 0.60601294, 0.65641482, 0.80975928
\\
\hline 
\end{tabular}
\label{tab:srj-schemes-c0}
\end{center}
\end{table}

\begin{table}[htbp!]
\caption{SRJ scheme parameters for $c = \frac{1}{10}$}
\begin{center}
\begin{tabular}{ |c|p{150mm}|  } 
\hline
$M$ & \textbf{SRJ scheme parameters} \\
\hline
2 & 0.56998789, 1.69859189 \\
\hline
3 & 0.53391192, 3.44601684, 0.92457397 \\
\hline
4 & 0.52000846, 5.84801488, 1.46164411, 0.70927878 \\
\hline
5 & 0.97045893, 8.85298484, 2.15794431, 0.51336698, 0.62598727 \\
\hline
6 & 3.00593438, 12.40210949,  0.78535133,  0.58484042,  1.30000478, 0.50970689 \\
\hline
7 & 16.43010567,  4.0003348, 1.69275624, 0.56135711, 0.98455299, 0.69414294, 0.50748255 \\
\hline
8 & 0.50603201, 20.86770462, 5.13612346, 0.83069204, 1.21911816, 2.14633476, 0.54663294, 0.64177848 \\
\hline
9 & 25.6439171, 6.40794264, 0.5050344, 0.74261526, 2.6591608, 0.53676794, 0.99056047, 1.48706252, 0.60864258 \\
\hline
10 & 30.68829417, 1.17192817, 1.78726047, 0.52982376, 0.50431918, 0.86044488, 7.8099309, 3.22986843, 0.68683046, 0.58621444 \\
\hline
11 & 0.50378925, 35.93393042, 0.99363246, 0.57027053, 0.64897316, 0.52474653, 3.85735463, 9.33630714, 0.7782492, 1.37385524, 2.11904043 \\
\hline
12 & 0.50338573, 41.31741418, 0.52091877, 0.62195236, 4.54033478, 0.55849849, 0.88139433, 1.1412681, 10.98051128, 0.72248676, 2.48180885, 1.59573262
\\
\hline
13 & 0.60191355, 0.68265354, 2.87506858, 5.27753675, 0.54954273, 0.50307142, 0.51796, 1.30280397, 12.73581301, 46.78097204, 0.80539615, 1.8371467, 0.99541991
\\
\hline
14 & 0.50282191, 52.2735576, 1.11983233, 14.595698, 0.58660159, 0.75115784, 0.65307138, 1.47791449, 0.89692809, 6.06783403, 0.51562534, 3.29845681, 2.09783992, 0.54256332
\\
\hline
15 & 0.51374938, 1.66627373, 6.90948562, 0.53701048, 0.82671796, 0.50262045, 0.63040987, 16.55218436, 57.74880892, 3.75132855, 0.99656637, 0.57460771, 1.25425115, 0.71085879, 2.37743854 \\
\hline
16 & 0.53251634, 1.10397129, 0.77459073, 0.50245545, 1.86764801, 63.16776026, 0.67996696, 0.56502152, 0.612616, 1.39841661, 0.51221886, 0.90885822, 7.80091195, 2.67565505, 18.59768396, 4.23314259 \\
\hline
17 & 4.74387164, 8.74143558, 0.59837029, 0.52882897, 0.50231873, 0.73467077, 0.55723606, 2.08203762, 0.6557056, 0.84386349, 0.99731915, 1.21899192, 0.51095483, 20.72672602, 68.50019668, 2.99252749, 1.55226259 \\
\hline
18 & 0.5022041, 0.55081535, 1.34139356, 0.58675992, 3.32758594, 73.7175464, 1.71554614, 1.09183698, 0.52576209, 2.30913802, 0.70328863, 0.63624059, 0.79408952, 22.93016201, 5.28260568, 9.72872304, 0.50989774, 0.91834936 \\
\hline
19 & 1.1922696, 1.47107281, 0.754984, 0.5090049, 5.84897369, 0.50210705, 0.85794803, 0.67810534, 1.88818571, 3.68068961, 10.76141436, 0.52318366, 25.20121758, 0.62035781, 0.99785036, 0.57716248, 78.79939956, 2.54886524, 0.54545533 \\
\hline
20 & 27.53285347, 83.72975055, 0.65754366, 11.83795003, 0.54093216, 0.81052691, 1.60793981, 6.44250844, 4.05164906, 0.50202419, 0.50824398, 0.72360616, 0.92605018, 0.5209947, 1.29850579, 1.08222029, 0.56913027, 2.8011136, 0.60720895, 2.0700968
\\
\hline 
\end{tabular}
\label{tab:srj-schemes-c1_10}
\end{center}
\end{table}

\begin{table}[htbp!]
\caption{SRJ scheme parameters for $c = \frac{1}{5}$}
\begin{center}
\begin{tabular}{ |c|p{150mm}| } 
\hline
$M$ & \textbf{SRJ scheme parameters} \\
\hline
2 & 1.67329868, 0.57289381 \\
\hline
3 & 3.30826695, 0.5373787, 0.92457405 \\
\hline
4 & 0.52365056, 0.71207692, 1.44990317, 5.42377641 \\
\hline
5 & 0.51708553, 0.62939828, 7.87621952, 2.11836778, 0.97045893 \\
\hline
6 & 10.52706852, 2.91384798, 1.29353361, 0.78773197, 0.58845872, 0.51346561 \\
\hline
7 & 0.51126504, 0.69728394, 3.82161829, 1.67437171, 0.56506573, 0.98455481, 13.25519709 \\
\hline
8 & 0.83271745, 0.64523131, 4.82675295, 0.55038743, 2.10859708, 1.21478179, 0.50982969, 15.96393814 \\
\hline
9 & 0.99054685, 0.50884231, 0.54054782, 0.61224365, 0.74545738, 1.47574166, 2.59243292, 18.58247992, 5.91397015 \\
\hline
10 & 0.50813454, 0.6900735, 1.1686967, 0.533621, 21.06405733, 0.58990017, 7.06842827, 1.76565972, 0.86219325, 3.1223753  \\
\hline
11 & 0.50761003, 0.65242796, 23.38109951, 0.57400444, 8.27539085, 0.52855467, 2.08301117, 0.99363078, 1.36587294,  3.69478206, 0.78083252 \\
\hline
12 & 1.13872598, 0.52473523, 2.42647806, 0.72552509, 0.50721069, 0.56226476, 9.5212926, 25.52113096, 0.62553501,  0.88294173, 1.58117851, 4.30608714 \\
\hline
13 & 1.81376857, 10.79255978, 27.48184901, 4.95240233, 2.79456256, 1.29671323, 0.6859496, 0.60557162, 0.80775937, 0.55332843, 0.52178179, 0.50689955, 0.9954343 \\
\hline
14 & 0.51945011, 0.50665244, 1.11768542, 0.54635948, 1.46703489, 3.18587451, 5.63006462, 2.06294125, 0.7539773, 12.07727937, 0.59030279, 0.89827566, 29.26841843, 0.65651632 \\
\hline
15 & 13.36645552, 30.89028993, 0.71398235, 0.50645307, 0.82887223, 0.51757784, 0.57834595, 0.63396522, 1.24931321,  0.54081793, 6.33555139, 1.64934075, 2.32817313, 0.99656232, 3.59938118 \\
\hline
16 & 0.91008368, 32.35924423, 0.56878365, 0.61624256, 0.51604983, 2.6086148, 1.38986873, 0.68328291, 1.10216036, 7.06508676, 0.77723577, 0.50628973, 4.0331308, 1.84305372, 14.64832029, 0.536331  \\
\hline
17 & 0.60204809, 33.68750902, 1.5391942, 0.56101623, 0.73763853, 0.50615441, 1.21487444, 0.51478793, 2.90385064,  0.65915403, 15.91870643, 0.99732293, 7.81500258, 4.48656154, 2.04803864, 0.53264945, 0.84585323 \\
\hline
18 & 17.16755739, 8.58203324, 0.50604092, 1.33437917, 0.70647254, 2.26371203, 0.91946802, 34.88850558, 3.21303678, 4.95751243, 0.51373245, 0.52958669, 0.63978006, 1.69691089, 0.59047345, 0.55460818, 1.09026941, 0.79657302 \\
\hline
19 & 0.99785451, 0.54925773, 18.39057107, 0.52701147, 9.36281884, 0.51284094, 2.48979141, 1.18873926, 5.4448419, 35.97445695, 3.53558757, 1.8628211, 1.4605054, 0.68144074, 0.62396284, 0.75780222, 0.85979157, 0.58090237, 0.50594485 \\
\hline
20 & 0.57288937, 36.95702075, 0.61086135, 0.50586281, 0.66098678, 3.87084352, 1.29256497, 2.03669706, 0.81285875, 0.72665698, 10.15421635, 19.58334469, 1.59307571, 1.08082848, 0.51208106, 5.94715006, 0.54474155, 2.72591983, 0.5248249, 0.92707131
\\
\hline
\end{tabular}
\label{tab:srj-schemes-c1_5}
\end{center}
\end{table}

\begin{table}[htbp!]
\caption{SRJ scheme parameters for $c = \frac{1}{3}$}
\begin{center}
\begin{tabular}{ |c|p{150mm}| } 
\hline
$M$ & \textbf{SRJ scheme parameters} \\
\hline
2 & 0.58009423, 1.61475671 \\
\hline
3 & 0.54601048, 3.01484954, 0.92457408 \\
\hline
4 & 0.53273474, 4.60959599, 0.71897996, 1.42207901 \\
\hline
5 & 2.02782143, 6.20847021, 0.52636835, 0.97045884, 0.63786073 \\
\hline
6 & 0.79358883, 1.27804404, 2.71159702, 0.52285298, 7.69468612, 0.59745813 \\
\hline
7 & 9.01421102, 0.70504972, 3.44669363, 1.63122788, 0.57430277, 0.5207143, 0.98455473 \\
\hline
8 & 0.83769031, 1.20435204, 0.65379547, 2.02203906, 0.55974784, 0.51931866, 4.20943149, 10.15428711 \\
\hline
9 & 0.5183584, 0.62120265, 0.54998076, 1.44897399, 0.99056019, 0.7524925, 11.12451034, 2.44332527, 4.97973358 \\
\hline
10 & 0.51766971, 5.74097156, 11.94430379, 0.86647514, 1.71513908, 1.16092303, 0.69810077, 2.88849809, 0.5430979, 0.59907044 \\
\hline
11 & 0.58330617, 1.34685534, 0.53806264, 0.66100147, 0.78720632, 6.48228504, 12.63494347,  0.51715946,  2.00027229,  3.35099357, 0.9936465 \\
\hline
12 & 7.19345185, 13.21763033, 2.301418, 0.51677077, 0.53426362, 0.73301243, 0.88670167, 0.6344228, 1.54674947, 3.82533345, 1.13253407, 0.57164568 \\
\hline
13 & 2.61594525, 0.51646796, 0.53132553, 0.69409679, 7.86844778, 0.61465857, 4.30651387, 0.81354555, 0.99544263, 13.71084518, 1.28204135, 0.56276371, 1.75922945 \\
\hline
14 & 14.12921075, 1.11258401, 1.44134806, 2.94184398, 0.55583452, 4.78908121, 0.51622752, 0.66508152, 0.90163402, 8.50561058, 0.76097013, 1.9828288, 0.59952601, 0.52900621 \\
\hline
15 & 14.48642564, 9.10177106, 0.55031875, 0.51603341, 0.58765634, 0.52714239, 0.83417426, 1.23747314, 0.72170963,  5.26980847, 0.99656684, 1.60955139, 0.64279542, 2.21639217, 3.27638036 \\
\hline
16 & 0.78377061, 0.62526217, 0.52562193, 0.54585339, 1.09781681, 1.78590467, 14.79284308, 0.57816083, 2.4586654, 0.91310944, 9.65718979, 1.36954054, 3.61745717, 0.51587446, 0.69150587, 5.74512382 \\
\hline
17 & 0.52436516, 10.17291044, 0.61119272, 0.51574268, 1.96968648, 6.21210459, 1.20494873, 0.99732608, 0.66770064, 15.05697411, 2.7084423, 0.74496189, 0.54218565, 3.96308173, 1.50827045, 0.57043683, 0.85073018 \\
\hline
18 & 0.80268015, 0.92219534, 4.31145962, 1.08646238, 0.53913484, 2.96453265, 0.71434257, 0.56406383, 0.52331433, 0.64856217, 15.28572255, 2.16021516, 6.66831474, 10.65068252, 1.31758874, 1.65317168, 0.5156322, 0.59971285 \\
\hline
19 & 0.51553868, 7.1121415, 0.76474781, 11.09221356, 4.66077265, 0.52242668, 0.63291549, 15.48494337, 1.8037996, 0.99785572, 2.35677891, 3.22588957, 1.18021839, 0.55874009, 0.53656909, 0.86430683, 1.43538507, 0.59021236, 0.68969701 \\
\hline 
20 & 0.61994079, 0.5542447, 0.92956967, 0.66952289, 0.58225288, 0.73419468, 0.52167002, 15.65937926, 3.49151906, 5.00938372, 1.27833205, 1.95972847, 0.51545881, 11.49957143, 7.54235807, 0.53439011, 1.07745477, 1.55800776, 2.5586843, 0.81859311 \\
\hline
\end{tabular}
\label{tab:srj-schemes-c1_3}
\end{center}
\end{table}

\begin{table}[htbp!]
\caption{SRJ scheme parameters for $c = \frac{1}{2}$}
\begin{center}
\begin{tabular}{ |c|p{150mm}| } 
\hline
$M$ & \textbf{SRJ scheme parameters} \\
\hline
2 & 0.59563557, 1.50541872 \\
\hline
3 & 0.56484541, 0.9245688, 2.54605612 \\
\hline
4 & 0.73375507, 3.51442902, 1.36762629, 0.55263795 \\
\hline
5 & 0.65617571, 0.54674458, 0.9704589, 4.31270689, 1.86254927\\
\hline
6 & 1.24704532, 2.37008646, 0.80603023, 0.61704862, 4.9370665, 0.54348006 \\
\hline
7 & 0.54149044, 0.59447669, 5.41652909, 1.54840168, 0.9845544, 0.72173662, 2.86349551 \\
\hline
8 & 0.54019056, 5.78430508, 3.32624599, 0.67231978, 0.58023364, 1.86326932, 0.84820779, 1.18325486 \\
\hline
9 & 6.06843828, 3.74978583, 0.5392957, 0.99056864, 2.18266908, 1.39632673, 0.64065916, 0.57065124, 0.76751184 \\
\hline
10 & 0.56388951, 6.29037877, 0.61905175, 4.13107823, 0.53865334, 2.49956211, 1.6191135, 0.7153856, 1.14510986, 0.87552324 \\
\hline
11 & 0.67952702, 6.46575637, 1.84779287, 0.53817747, 0.99364015, 0.8007591, 2.80842724, 0.55893338, 0.60360767, 1.30897648, 4.47091216 \\
\hline
12 & 6.60623628, 3.10568012, 0.59215413, 0.89462012, 0.74905399, 0.65370652, 4.77144259, 0.53781461, 1.11986881, 2.07889593, 0.55519275, 1.47989517 \\
\hline
13 & 2.30978735, 1.25268301, 0.58341686, 0.82583806, 3.38846892, 1.65574798, 5.03656248, 0.99544475, 0.63442357, 6.71996498, 0.53753215, 0.71164398, 0.55229695 \\
\hline
14 & 0.5500114, 6.81350362, 0.77590957, 0.53730728, 1.10217112, 0.57658532, 1.83507103, 1.39055807, 2.53760691, 0.68356729, 5.26917182, 3.65606827, 0.90860206, 0.61962197 \\
\hline
15 & 5.47370042, 6.89103876, 0.84546621, 0.54816904, 2.76076389, 0.53712712, 0.66195405, 0.99643676, 0.73825278, 0.60796165, 1.21381076, 3.90725337, 1.53236617, 0.57115304, 2.01615777 \\
\hline
16 & 0.56673831, 0.6448415, 0.53697797, 0.7976593, 5.6558776, 4.13962774, 1.6777469, 0.91941575, 2.97955895, 6.95497755, 0.59864266, 1.08883637, 0.54667133, 1.32911154, 2.19646221, 0.70920553 \\
\hline
17 & 0.86073193, 5.80924455, 0.9978438, 3.18482579, 2.3803012, 7.01207131, 0.53685725, 4.36359145, 0.76081671, 1.44933689, 0.54542277, 1.82234187, 0.56313151, 0.63115087, 0.59100793, 1.18398732, 0.6860663 \\
\hline
18 & 4.56120355, 0.56009249, 5.95395338, 0.61982757, 0.8157569, 1.97175496, 3.3908333, 1.07890479, 0.58477151, 0.73119493, 0.92774785, 2.55624736, 1.56992948, 1.28360488, 7.05683188, 0.53675096, 0.54439601, 0.6676193 \\
\hline
19 & 4.74590573, 0.53666476, 2.12088425, 0.54351562, 3.58685822, 0.57950529, 6.08052385, 0.70750994, 0.87388224, 0.77954849, 1.69284254, 0.65233427, 7.09590848, 0.9977728, 2.73045577, 1.38567105, 0.55756252, 0.61053146, 1.16304593 \\
\hline 
20 &  2.90348415, 0.63969696, 0.93492526, 1.25009592, 7.13078781, 4.92003011, 0.830678, 0.60269106, 0.6879351, 6.19014037, 1.07030842, 1.81807115, 3.7710704, 0.55540038, 0.53658931, 0.7504016, 1.48936404, 0.54277116, 0.5750778, 2.26815286 \\
\hline
\end{tabular}
\label{tab:srj-schemes-c1_2}
\end{center}
\end{table}

\end{document}